\documentclass[11pt]{amsart}
\usepackage[english]{babel}
\usepackage{amsthm} 
\usepackage{amssymb}
\usepackage{amsmath}
\usepackage{hyperref}
\usepackage{amsfonts}
\usepackage{verbatim}
\usepackage{color}
\usepackage{stmaryrd}
\usepackage{fancyhdr}
\usepackage{mathrsfs}
\usepackage{graphicx}
\usepackage{pdfsync}

\usepackage{esint} 

\newtheorem{theorem}{Theorem}[section]
\newtheorem{corollary}[theorem]{Corollary}
\newtheorem{definition}[theorem]{Definition}
\newtheorem{lemma}[theorem]{Lemma}
\newtheorem{proposition}[theorem]{Proposition}
\newtheorem{remark}{Remark}
\numberwithin{equation}{section}

\newcommand{\rr}{\mathbb{R}}
\newcommand{\cc}{\mathbb{C}}

\newcommand{\eps}{\varepsilon}

\def\un{{\mathrm{1~\hspace{-1.4ex}l}}}

\def\un{{\mathrm{1~\hspace{-1.4ex}l}}}

\usepackage{varioref}
 
\def\Rom#1{\uppercase\expandafter{\romannumeral #1}}

\def\polhk#1{\setbox0=\hbox{#1}{\ooalign{\hidewidth \lower1.5ex\hbox{`}\hidewidth\crcr\unhbox0}}}

\author{Karine \textsc{Beauchard}, Karel \textsc{Pravda-Starov}}

\address{\noindent \textsc{Karine Beauchard, IRMAR, CNRS UMR 6625, \'Ecole normale sup\'erieure de Rennes, Avenue Schumann, 35170 Bruz, France}}
\email{karine.beauchard@ens-rennes.fr}
\address{\noindent \textsc{Karel Pravda-Starov, IRMAR, CNRS UMR 6625, Universit\'e de Rennes 1, Campus de Beaulieu, 263 avenue du G\'en\'eral Leclerc, CS 74205,
35042 Rennes cedex, France}}
\email{karel.pravda-starov@univ-rennes1.fr}


\title{Null-controllability of non-autonomous Ornstein-Uhlenbeck equations}
\date{}

\begin{document}

\begin{abstract}
We study the null-controllability of parabolic equations associated to non-autonomous Ornstein-Uhlenbeck operators. When a Kalman type condition holds for some positive time $T>0$, these parabolic equations are shown to enjoy a Gevrey regularizing effect at time $T>0$. Thanks to this regularizing effect, we prove by adapting the Lebeau-Robbiano method that these parabolic equations are null-controllable in time $T>0$ from control regions, for which null-controllability is classically known to hold in the case of the heat equation. 
\end{abstract}

\maketitle


\section{Introduction}

\subsection{Null-controllability of degenerate parabolic equations}
We aim in this work at studying the null-controllability of parabolic equations posed on the whole space $\rr^n$ and controlled by a source term locally distributed on an open subset $\omega \subset \mathbb{R}^n$,
\begin{equation}\label{syst_general}
\left\lbrace \begin{array}{l}
\partial_tf(t,x) -  \frac{1}{2}\textrm{Tr}\big(A(t)A(t)^T\nabla_x^2f(t,x)\big) - \big\langle B(t)x, \nabla_x f(t,x)\big\rangle=u(t,x)\un_{\omega}(x),  \\
f|_{t=0}=f_0 \in L^2(\rr^n),                                       
\end{array}\right.
\end{equation}
associated to non-autonomous Ornstein-Uhlenbeck operators
$$P(t)=\frac{1}{2}\textrm{Tr}\big(A(t)A(t)^T\nabla_x^2\big) + \big\langle B(t)x, \nabla_x \big\rangle, \quad   x \in \mathbb{R}^n,$$
where $A, B \in C^{\infty}(I,M_{n}(\rr)),$ 
are smooth mappings with values in real $n \times n$ matrices, with $I$ being an open interval of $\rr$ containing zero, $A(t)^T$ standing for the transpose matrix of $A(t)$ and $\un_{\omega}$ denoting the characteristic function of the set $\omega$.
Here, the notation $\langle A,B \rangle$ stands for the scalar operator
\begin{equation}\label{not}
\langle A,B\rangle=\sum_{j=1}^nA_jB_j,
\end{equation}
when $A=(A_1,...,A_n)$ and $B=(B_1,...,B_n)$ are vector-valued operators. Notice that $\langle A,B\rangle\neq\langle B,A\rangle$ in general, since e.g., 
$\langle \nabla_x,Bx\rangle=\langle Bx,\nabla_x\rangle +\textrm{Tr}(B)$, as for any $u \in \mathscr{S}'(\rr^n)$,
$$\langle \nabla_x,Bx\rangle u=\sum_{j=1}^n \partial_{x_j}\big((Bx)_ju\big)=\sum_{j=1}^n (Bx)_j\partial_{x_j}u+\sum_{j=1}^n \partial_{x_j}\big((Bx)_j\big)u=\langle Bx,\nabla_x\rangle u +\textrm{Tr}(B)u.$$
The null-controllability of the parabolic equation \eqref{syst_general} is defined as follows:

\medskip

\begin{definition} [Null-controllability] Let $T>0$ and $\omega$ be an open subset of $\mathbb{R}^n$. 
Equation \eqref{syst_general} is said to be {\em null-controllable from the set $\omega$ in time} $T$ if, for any initial datum $f_0 \in L^{2}(\mathbb{R}^n)$, there exists $u \in L^2((0,T)\times\mathbb{R}^n)$, supported in $[0,T]\times\omega$, such that the mild solution of \eqref{syst_general} satisfies $f(T,\cdot)=0$.
\end{definition}

\medskip

By the Hilbert uniqueness method (see Proposition~\ref{Prop:HUM} in appendix), the null controllability of the equation \eqref{syst_general} in time $T>0$ is equivalent to the observability of the adjoint system in time $T>0$, 
\begin{equation} \label{adj_general}
\left\lbrace \begin{array}{ll}
\partial_tg(t,x) - P(T-t)^*g(t,x)=0\,, \quad & x \in \mathbb{R}^n\,, \\
g|_{t=0}=g_0 \in L^2(\rr^n).
\end{array}\right.
\end{equation}
We recall that the notion of observability is defined as follows:

\medskip

\begin{definition} [Observability] Let $T>0$ and $\omega$ be an open subset of $\mathbb{R}^n$. 
Equation \eqref{adj_general} is said to be {\em observable in the set $\omega$ in time} $T$ if there exists a constant $C>0$ such that,
for any initial datum $g_0 \in L^{2}(\mathbb{R}^n)$, the mild solution of \eqref{adj_general} satisfies
\begin{equation}\label{eq:observability}
\int\limits_{\mathbb{R}^n} |g(T,x)|^{2} dx  \leqslant C \int\limits_{0}^{T} \Big(\int\limits_{\omega} |g(t,x)|^{2} dx\Big) dt\,.
\end{equation}
\end{definition}

\medskip

An important open problem at the core of current investigations is to understand to which extent the null-controllability or observability results known for uniformly parabolic equations still hold for degenerate parabolic equations of hypoelliptic type.

For equations posed on bounded domains, some progress have been made. In the case of the heat equation on a bounded domain $\Omega$ with Dirichlet boundary conditions, it is well-known that observability holds in arbitrary positive time $T>0$, with any non-empty open set $\omega$, see~\cite[Theorem~3.3]{Fattorini-Russel}, \cite{Fursikov-Imanuvilov-186} and \cite{Lebeau-Robbiano}. Degenerate parabolic equations exhibit a wider range of behaviours. Indeed, observability may hold true, or not, depending on the strength of the degeneracy. This feature is well understood for parabolic equations that degenerate on the domain boundary, see~\cite{Ala-Can-Fra,Can-Fra-Roc_2,Can-Fra-Roc_1, Cannarsa-V-M-ADE, Cannarsa-V-M-SIAM, Martinez-Vancost-JEE-2006} in the one-dimensional case, and~\cite{Cannarsa-V-M-CRAS} for the multi-dimensional one. Furthermore, a positive minimal time may be required to get observability, see the works~\cite{Grushin, MR3162108} in the case of the Grushin equation, \cite{KB_Can_Heisenberg} for the Heisenberg heat equation, and~\cite{MR3163490} for the Kolmogorov equation. This minimal time is actually related to localization properties of eigenfunctions. Finally, a geometric control condition may also be required for the observability inequality to hold~\cite{KB_Helffer}.

On the other hand, the understanding of the null-controllability and observability for degenerate parabolic equations of hypoelliptic type posed on the whole space is still at an earlier stage. For the heat equation on the whole space
\begin{equation}\label{heat_ws}
\left( \partial_t - \Delta_x \right)f(t,x)=u(t,x)\un_{\omega}(x)\,, \quad (t,x)\in(0,T)\times\mathbb{R}^n,
\end{equation}
no necessary and sufficient condition on the control region $\omega$ is known for null-controllability to hold in any positive time. The condition
$$\sup_{x \in \mathbb{R}^n}d(x,\omega)<+\infty,$$
is shown in~\cite[Theorem 1.11]{Miller_BSM2005} to be necessary for null-controllability to hold in any positive time. On the other hand, the following sufficient condition  
\begin{equation}\label{hyp_omega}
\exists \delta, r >0, \forall y \in \mathbb{R}^n\,, \exists y' \in \omega,\quad B(y',r) \subset \omega \text{ and } |y-y'|<\delta\,,
\end{equation}
is given in~\cite{Miller_unbounded} for null-controllability to hold from the open set $\omega \subset \rr^n$ in any positive time. 
The very same condition is shown in~\cite{LeRousseau_Moyano} to be sufficient, when the observability set $\omega=\omega_x \times \omega_v$ has a Cartesian structure with $\omega_x$ and $\omega_v$ open sets both satisfying (\ref{hyp_omega}), for the null-controllability of the Kolmogorov equation
\begin{equation}\label{koleq}
\left\lbrace \begin{array}{ll}
(\partial_t + v \cdot \nabla_x-\Delta_v)f(t,x,v)=u(t,x,v)\un_{\omega}(x,v)\,, \quad &  (x,v) \in \mathbb{R}^{n_1} \times \mathbb{R}^{n_2},\\
f|_{t=0}=f_0 \in L^2(\rr_{x}^{n_1} \times \rr_{v}^{n_2}),                                       &  
\end{array}\right.
\end{equation}
to hold in any positive time. This result of null-controllability of the Kolmogorov equation was then extended by Zhang~\cite{zhang} without the constraint on the Cartesian structure of the observability set.

More generally, the condition (\ref{hyp_omega}) was next shown to be sufficient for the null controllability of all hypoelliptic Ornstein-Uhlenbeck equations~\cite{KarKar1} (Theorem~1.3),
\begin{equation} \label{syst_LR}
\left\lbrace \begin{array}{ll}
\partial_t f(t,x) - \frac{1}{2}\textrm{Tr}\big(Q\nabla_x^2 f(t,x)\big) - \big\langle Bx, \nabla_x f(t,x)\big\rangle=u(t,x)\un_{\omega}(x)\,, \quad &  x \in \mathbb{R}^n,\\
f|_{t=0}=f_0 \in L^2(\rr_x^n),                                       &  
\end{array}\right.
\end{equation}
where $Q$, $B$ are real $n \times n$-matrices, with $Q$ symmetric positive semidefinite, satisfying the Kalman rank condition
\begin{equation}\label{kal1}
\textrm{Rank}[B|Q^{\frac{1}{2}}]=n,
\end{equation} 
where 
$$[B|Q^{\frac{1}{2}}]=[Q^{\frac{1}{2}},BQ^{\frac{1}{2}},\dots, B^{n-1}Q^{\frac{1}{2}}],$$ 
is the $n\times n^2$ matrix obtained by writing consecutively the columns of the matrices $B^jQ^{\frac{1}{2}}$, with $Q^{\frac{1}{2}}$ the symmetric positive semidefinite matrix given by the square root of $Q$. This general result allows in particular while taking $Q=2I_n$ and $B=0$, to recover the result of null-controllability of the heat equation (\ref{heat_ws}), and by taking 
$$Q=\left(\begin{array}{ll}
0 & 0 \\ 0 & 2I_d
\end{array}\right), \qquad
B=\left(\begin{array}{ll}
0 & -I_d \\ 0 & 0
\end{array}\right), \qquad n=2d,$$
to recover the result of null-controllability of the Kolmogorov equation (\ref{koleq}) established in~\cite{zhang}.

In the present work, we aim at extending this result of null-controllability of Ornstein-Uhlenbeck equations in the non-autonomous case
$$\left\lbrace \begin{array}{l}
\partial_tf(t,x) -  \frac{1}{2}\textrm{Tr}\big(A(t)A(t)^T\nabla_x^2f(t,x)\big) - \big\langle B(t)x, \nabla_x f(t,x)\big\rangle=u(t,x)\un_{\omega}(x), \quad x \in \rr^n, \\
f|_{t=0}=f_0 \in L^2(\rr^n),                                       
\end{array}\right.$$
with $A, B \in C^{\infty}(I,M_{n}(\rr))$, when a generalization of the Kalman condition to the time-varying case holds.

\bigskip

\subsection{Statement of the main result}

Let $I$ be an open interval of $\rr$ containing zero and $P(t)$ be the non-autonomous Ornstein-Uhlenbeck operator
\begin{multline}\label{NAOU}
P(t)=\frac{1}{2}\textrm{Tr}\big(A(t)A(t)^T\nabla_x^2\big) + \big\langle B(t)x, \nabla_x \big\rangle\\
=\frac{1}{2}\sum_{i,j,k=1}^{n}a_{i,k}(t)a_{j,k}(t)\partial_{x_i,x_j}^2+\sum_{i,j=1}^nb_{i,j}(t)x_j\partial_{x_i}, \quad   x \in \mathbb{R}^n,
\end{multline}
where $$A=(a_{i,j})_{1 \leq i,j \leq n},\quad B=(b_{i,j})_{1 \leq i,j \leq n} \in C^{\infty}(I,M_{n}(\rr)),$$ 
are smooth mappings with values in real $n \times n$ matrices. We define by induction the sequence of smooth mappings $(\tilde{A}_k)_{k \geq 0} \in C^{\infty}(I,M_{n}(\rr))^{\mathbb{N}}$ as
\begin{equation}\label{asd2}
\forall t \in I, \quad \tilde{A}_0(t)=A(t),
\end{equation}
\begin{equation}\label{asd3}
\forall k \geq 0, \forall t \in I, \quad \tilde{A}_{k+1}(t)= \frac{d}{dt}\tilde{A}_{k}(t)+B(t)\tilde{A}_{k}(t).
\end{equation}
We consider the following generalization of the Kalman condition to the time-varying case
\begin{equation}\label{Kalman_time}
\exists t_0 \in I, \quad \textrm{Span}\{\tilde{A}_{k}(t_0)x : \ x \in \rr^n,\ k \geq 0\}=\rr^n.
\end{equation}
This condition (\ref{Kalman_time}) was shown by Chang~\cite{chang} and by Silverman and Meadows~\cite{silverman} to be sufficient for the controllability of the linear control system $\dot{x}=-B(t)x+A(t)u$ on the interval $I$. As noticed in~\cite[p. 11]{coron_book}, it is important to notice that the two following vector spaces 
$$\textrm{Span}\{\tilde{A}_{k}(t_0)x : \ x \in \rr^n,\ k \geq 0\} \neq \textrm{Span}\{\tilde{A}_{k}(t_0)x : \ x \in \rr^n,\ 0 \leq k \leq n-1\},$$
are in general distinct, as contrary to the constant case, the Cayley-Hamilton theorem does not apply. However, it was proved by Coron in~\cite{coron_book} (Proposition~1.19) that when the condition (\ref{Kalman_time}) holds at some time $t_0 \in I$, then there exists a positive constant $\eps>0$ such that 
\begin{equation}\label{Kalman_time11}
\forall t \in I \cap (t_0-\eps,t_0+\eps)\setminus \{t_0\}, \quad \textrm{Span}\{\tilde{A}_{k}(t)x : \ x \in \rr^n,\ 0 \leq k \leq n-1\}=\rr^n.
\end{equation}
This assertion (\ref{Kalman_time11}) can be reformulated as
\begin{equation}\label{Kalman_time2}
\forall t \in I \cap (t_0-\eps,t_0+\eps)\setminus \{t_0\}, \quad \textrm{Rank}[\tilde{A}_{0}(t),\tilde{A}_{1}(t),\dots,\tilde{A}_{n-1}(t)]=n.
\end{equation}
 This reduction to finitely many matrices in (\ref{Kalman_time11}) and (\ref{Kalman_time2}) is used in the proof of the following result of null-controllability of non-autonomous Ornstein-Uhlenbeck equations, which is the main result contained in this work:

\medskip

\begin{theorem} \label{meta_thm}
Let $I$ be an open interval of $\rr$ containing zero and $\omega$ be an open subset of $\mathbb{R}^n$ satisfying
$$\exists \delta, r >0, \forall y \in \mathbb{R}^n\,, \exists y' \in \omega,\quad B(y',r) \subset \omega \text{ and } |y-y'|<\delta.$$
When the Kalman 
type condition 
\begin{equation}\label{wer1}
\emph{\textrm{Span}}\{\tilde{A}_{k}(T)x : \ x \in \rr^n,\ k \geq 0\}=\rr^n,
\end{equation}
holds for some positive time $T>0$ belonging to the interval $I$, the non-autonomous Ornstein-Uhlenbeck equation posed in the $L^2(\rr^n)$-space
\begin{equation}\label{syst_LR1}
\left\lbrace \begin{array}{l}
\partial_t f(t,x) - \frac{1}{2}\emph{\textrm{Tr}}\big(A(t)A(t)^T\nabla_x^2 f(t,x)\big) - \big\langle B(t)x, \nabla_x f(t,x)\big\rangle=u(t,x)\un_{\omega}(x), \\
f|_{t=0}=f_0 \in L^2(\rr_x^n),                                       
\end{array}\right.
\end{equation}
is null-controllable from the set $\omega$ in time greater or equal to $T$.
\end{theorem}

\medskip

This result extends the result of null-controllability of hypoelliptic Ornstein-Uhlenbeck equations obtained in~\cite{KarKar1}. Indeed, we recover the result of null-controllability of the Ornstein-Uhlenbeck equation
$$\left\lbrace \begin{array}{ll}
\partial_t f(t,x) - \frac{1}{2}\textrm{Tr}\big(Q\nabla_x^2 f(t,x)\big) - \big\langle Bx, \nabla_x f(t,x)\big\rangle=u(t,x)\un_{\omega}(x)\,, \quad &  x \in \mathbb{R}^n,\\
f|_{t=0}=f_0 \in L^2(\rr_x^n),                                       &  
\end{array}\right.$$
where $Q$, $B$ are real $n \times n$-matrices, with $Q$ symmetric positive semidefinite, when taking $A(t)=Q^{\frac{1}{2}}$ and $B(t)=B$, since in this case the Kalman type condition (\ref{wer1}) reads as 
$$\textrm{Span}\{B^kQ^{\frac{1}{2}}x : \ x \in \rr^n,\ k \geq 0\}=\rr^n,$$
that is,
$$\textrm{Span}\{B^kQ^{\frac{1}{2}}x : \ x \in \rr^n,\ 0 \leq k \leq n-1\}=\rr^n,$$
thanks to the Cayley-Hamilton theorem, the latter condition being equivalent to the Kalman condition (\ref{kal1}).

The proof of Theorem~\ref{meta_thm} in this work is an adaptation of the proof given in~\cite{KarKar1} (Theorem~1.3). This proof follows the Lebeau-Robbiano strategy for establishing the observability inequality of the adjoint system (\ref{adj_general}) and relies on the fact that this Cauchy problem can be solved explicitly, see the appendix in Section~\ref{appendix}. 
Nevertheless, the Lebeau-Robbiano method cannot be directly applied in its usual form. This was already noticed in~\cite{KarKar1} in the constant case since Ornstein-Uhlenbeck semigroups do not commute with the Fourier frequency cutoff projections. This non-commutation accounts that even if some low frequencies could be steered to zero at some time, any passive control phase in the Lebeau-Robbiano method makes them all revive again. To overcome this lack of commutation, we take a key advantage of the Gevrey smoothing effect (see Corollary~\ref{corollary}) enjoyed by the solutions to the Cauchy problem (\ref{adj_general}), when the Kalman type condition (\ref{wer1}) holds.

This Kalman type condition (\ref{wer1}) is known~\cite{chang,silverman} to be sufficient for the controllability of the linear system 
\begin{equation}\label{ghj}
\dot{x}=-B(t)x+A(t)u,
\end{equation}
on any interval $[T-\eps,T]$, with $0<\eps \leq T$. In general, this condition is not necessary for the controllability of the system (\ref{ghj}) unless $n=1$, or when the matrices $A$ and $B$ are analytic in the time variable, see for instance the example on p.11-12 together with Exercise 1.23 on p.19 in~\cite{coron_book}. On the other hand, a standard necessary and sufficient condition for the controllability of (\ref{ghj}) is given by the invertibility of the associated controllability Gramian~\cite{coron_book} (Theorem~1.11), which is equivalent to the positive definiteness of the quadratic form 
$$\xi \in \rr^n \mapsto \int_0^{\eps}|A(T-s)^TR(0,s)^T\xi|^2ds,$$
related to the quadratic form appearing in Lemma~\ref{alg}. This invertibility of the controllability Gramian
is not sufficient for our proof of the null-controllability of the non-autonomous Ornstein-Uhlenbeck equation (\ref{syst_LR1}). Indeed, a key ingredient in our proof is the explicit time estimate from below of the quadratic form given by Lemma~\ref{alg}, which does not necessarily hold when the Gramian is invertible unless $n=1$, or when the matrices $A$ and $B$ are analytic in the time variable, since in these cases, the invertibility of the controllability Gramian is equivalent to (\ref{wer1}).

When the matrices $A$ and $B$ are solely $C^{\infty}$, the Kalman type condition (\ref{wer1}) sufficient for the null-controllability of the non-autonomous Ornstein-Uhlenbeck equation (\ref{syst_LR1}) is not in general necessary. We refer the reader to Section~\ref{counter} for a counter-example of a null-controllable non-autonomous Ornstein-Uhlenbeck equation for which the Kalman type condition (\ref{wer1}) does not hold. It would be most interesting to understand further if the Kalman type condition (\ref{wer1}) is necessary for the null-controllability of the non-autonomous Ornstein-Uhlenbeck equation (\ref{syst_LR1}) when the matrices $A$ and $B$ are supposed to be analytic in the time variable as done in the work~\cite{burgos} on time-dependent coupled linear parabolic systems.

\begin{remark} When the Kalman type condition (\ref{wer1}) is strengthened as
\begin{equation}\label{ghj2}
\exists T >0, T \in I, \exists m \geq 0, \quad \emph{\textrm{Span}}\{\tilde{A}_{k}(T)x : \ x \in \rr^n,\ 0 \leq k \leq m\}=\rr^n,
\end{equation}
Theorem~\ref{meta_thm} holds under the weaker regularity assumptions $A \in C^{2m+1}(I,M_n(\rr))$ and $B \in C^{2m}(I,M_n(\rr))$. Indeed, these limited regularity assumptions are sufficient to get the result of Lemma~\ref{alg} when condition (\ref{ghj2}) holds. The other parts of the proof do not use any further regularity assumptions.
\end{remark}

\subsection{Outline of the work}
The next section is devoted to the proof of Theorem~\ref{meta_thm}, whereas the appendix (Section~\ref{appendix}) is dedicated to well-posedness results of the homogeneous and inhomogeneous Cauchy problems (\ref{syst_general}) and (\ref{adj_general}) together with some recalls about the Hilbert uniqueness method in the non-autonomous case and an explicit counter-example to the necessity of the Kalman type condition for null-controllability of non-autonomous Ornstein-Uhlenbeck equations  in the $C^{\infty}$ setting.

\section{Proof of Theorem~\ref{meta_thm}}\label{sec:LR}

This section is devoted to the proof of Theorem~\ref{meta_thm}. We begin by noticing that the result of Theorem~\ref{meta_thm} is equivalent to the null-controllability of the equation 
\begin{equation} \label{syst_LR1.z}
\left\lbrace \begin{array}{l}
\partial_t f - \frac{1}{2}\textrm{Tr}\big(A(t)A(t)^T\nabla_x^2 f\big) - \big\langle B(t)x, \nabla_x f\big\rangle-\frac{1}{2}\textrm{Tr}\big(B(t)\big)f=u(t,x)\un_{\omega}(x)\,, \\
f|_{t=0}=f_0 \in L^2(\rr_x^n),                                        
\end{array}\right.
\end{equation}
from the set $\omega$ in time $T>0$, where $\omega$ is an open subset of $\mathbb{R}^n$ satisfying (\ref{hyp_omega}). 
We next observe that the $L^2$-adjoint of the operator 
$$\frac{1}{2}\textrm{Tr}\big(A(t)A(t)^T\nabla_x^2\big)+\big\langle B(t)x, \nabla_x\big\rangle+\frac{1}{2}\textrm{Tr}\big(B(t)\big),$$ 
is given by
\begin{multline*}
\Big(\frac{1}{2}\textrm{Tr}\big(A(t)A(t)^T\nabla_x^2\big)+\big\langle B(t)x, \nabla_x\big\rangle+\frac{1}{2}\textrm{Tr}\big(B(t)\big)\Big)^*\\
=\frac{1}{2}\textrm{Tr}\big(A(t)A(t)^T\nabla_x^2\big)-\big\langle B(t)x, \nabla_x\big\rangle-\frac{1}{2}\textrm{Tr}\big(B(t)\big).
\end{multline*}
By using the Hilbert uniqueness method (see Proposition~\ref{Prop:HUM} in appendix), the null-controllability of the equation (\ref{syst_LR1.z}) from the set $\omega$ in time $T>0$ is equivalent to the following observability estimate
\begin{equation}\label{dfg2}
\exists C>0, \forall g_0 \in L^2(\rr^n), \quad \|g(T)\|_{L^2(\rr^n)}^2 \leq C\int_{0}^T\|g(t)\|_{L^2(\omega)}^2dt,
\end{equation}
where $g$ is the mild solution to the Cauchy problem
\begin{equation} \label{adj_general1}
\left\lbrace \begin{array}{l}
\partial_tg(t,x) - \tilde{P}(t)g(t,x)=0, \quad x \in \rr^n,\\
g|_{t=0}=g_0 \in L^2(\rr^n).
\end{array}\right.
\end{equation}
with 
$$\tilde{P}(t)=\frac{1}{2}\textrm{Tr}\big(A(T-t)A(T-t)^T\nabla_x^2\big)-\big\langle B(T-t)x, \nabla_x\big\rangle-\frac{1}{2}\textrm{Tr}\big(B(T-t)\big).$$
In the following, we denote $R$ the resolvent of the time-varying linear system 
$$\overset{.}{X}(t)=B(T-t)X(t),$$ 
defined as the mapping
\begin{equation}\label{resol}
\begin{array}{ll}
R : [0,T] \times [0,T] & \rightarrow M_n(\rr), \\
\ \ \ \  (t_1,t_2)  & \mapsto R(t_1,t_2),
\end{array}
\end{equation}
such that for all $t_2 \in [0,T]$, the mapping $t_1 \in [0,T] \mapsto R(t_1,t_2) \in M_n(\rr)$ is the solution of the Cauchy problem
$$\left\lbrace \begin{array}{l}
\overset{.}{X}(t)=B(T-t)X(t), \\
X(t_2)=I_n,
\end{array}\right.
$$
where $I_n$ denotes the identity matrix.

Next lemma is instrumental in the proof of Theorem~\ref{meta_thm} and allows to deduce the Gevrey smoothing effect enjoyed by the mild solutions to the Cauchy problem (\ref{adj_general1}), when the Kalman type condition (\ref{wer1}) holds:

\medskip

\begin{lemma}\label{alg}
Let $I$ be an open interval of $\rr$ containing zero and $A, B \in C^{\infty}(I,M_{n}(\rr))$. When the Kalman 
type condition (\ref{Kalman_time}) holds for some positive time $T>0$ belonging to the interval $I$,
then there exists a positive constant $0<\eps <T$ such that 
$$\forall T-\eps \leq t \leq T-\frac{\eps}{2}, \quad \emph{\textrm{Span}}\{\tilde{A}_{k}(t)x : \ x \in \rr^n,\ 0 \leq k \leq n-1\}=\rr^n.$$
Furthermore, there exist positive constants $c>0$, $0< \tilde{\eps} \leq \frac{\eps}{2}$ such that 
$$\forall \frac{\eps}{2} \leq t \leq \tau \leq \frac{\eps}{2}+\tilde{\eps}, \forall \xi \in \rr^{n}, \quad \int_{t}^{\tau}|A(T-s)^TR(t,s)^T\xi|^2ds \geq c(\tau-t)^{2n-1}|\xi|^2,$$
with $|\cdot|$ being the Euclidean norm on $\rr^{n}$.
\end{lemma}

\medskip 

\begin{proof} We first deduce from the Kalman type condition holding at time $T>0$, and~\cite{coron_book} (Proposition~1.19) that there exists a positive constant $0<\eps< \min(1,T)$ such that for all $T-\eps \leq t \leq T-\frac{\eps}{2}$,
\begin{equation}\label{Kalman_time1}
\textrm{Span}\{\tilde{A}_{k}(t)x : \ x \in \rr^n,\ 0 \leq k \leq n-1\}=\rr^n.
\end{equation}
We recall for instance from~\cite{coron_book} (Proposition~1.5) that the resolvent satisfies the following properties 
\begin{equation}\label{asd1}
\forall t, \tau \in [0,T], \ R(t,\tau)R(\tau,t)=I_n, \quad \forall t,\tau \in [0,T], \ (\partial_2 R)(t,\tau)=-R(t,\tau)B(T-\tau).
\end{equation}
We notice from (\ref{asd2}), (\ref{asd3}) and (\ref{asd1}) that
\begin{equation}\label{asd4}
\forall k \geq 0, \forall t, \tau \in [0,T], \quad \frac{d^{k}}{d\tau^{k}}\big(A(T-\tau)^TR(t,\tau)^T\big)=(-1)^k\tilde{A}_k(T-\tau)^TR(t,\tau)^T.
\end{equation}
We consider the function
$$f_{\xi}(t,\tau)=\int_{t}^{\tau}|A(T-s)^TR(t,s)^T\xi|^2ds, \quad t, \tau \in [0,T],$$
depending on the parameter $\xi \in \rr^{n}$.
According to (\ref{asd4}), we easily check by the Leibniz formula that 
\begin{multline}\label{re1}
\forall n \geq 0, \forall t, \tau \in [0,T], \\ (\partial_2^{n+1}f_{\xi})(t,\tau)=(-1)^n\sum_{k=0}^n\binom{n}{k} \langle \tilde{A}_{n-k}(T-\tau)^TR(t,\tau)^T\xi,
\tilde{A}_k(T-\tau)^TR(t,\tau)^T\xi\rangle,
\end{multline}
where $\langle \cdot,\cdot \rangle$ denotes the Euclidean scalar product on $\rr^{n}$. 
We deduce from (\ref{Kalman_time1}) that 
\begin{equation}\label{asd5}
\forall  \frac{\eps}{2} \leq t \leq \eps, \quad \textrm{Rank}\big[\tilde{A}_{0}(T-t),\tilde{A}_{1}(T-t),\dots, \tilde{A}_{n-1}(T-t)\big]=n.
\end{equation}
We therefore have 
$$\forall \frac{\eps}{2} \leq t \leq \eps, \quad  \textrm{Ran}\big(\tilde{A}_{0}(T-t)\big)+\textrm{Ran}\big(\tilde{A}_{1}(T-t)\big)+...+\textrm{Ran}\big(\tilde{A}_{n-1}(T-t)\big)=\rr^{n}.$$
This implies that 
\begin{equation}\label{re2}
\forall \frac{\eps}{2} \leq t \leq \eps, \quad  \bigcap_{j=0}^{n-1}\textrm{Ker}\big(\tilde{A}_{j}(T-t)^T\big) \cap \rr^{n}=\{0\}.
\end{equation}
By induction, we easily check from (\ref{re1}) that for all $k \geq 0$, $\frac{\eps}{2} \leq t  \leq \eps$,
\begin{equation}\label{re3}
\forall 0 \leq l \leq 2k+1, \ (\partial_2^lf_{\xi})(t,\tau)|_{\tau=t}=0 \Longleftrightarrow  \xi \in \bigcap_{j=0}^{k}\textrm{Ker}\big(\tilde{A}_{j}(T-t)^T\big) \cap \rr^{n}.
\end{equation}
According to (\ref{re1}), (\ref{re2}) and (\ref{re3}), it follows that for all $\xi \in \rr^{n} \setminus \{0\}$, there exists $0 \leq k_{\xi} \leq n-1$ such that 
\begin{equation}\label{re0}
\forall 0 \leq j \leq 2k_{\xi}, \quad (\partial_2^jf_{\xi})\Big(\frac{\eps}{2},\frac{\eps}{2}\Big)=0
\end{equation}
and
\begin{equation}\label{re0b}
(\partial_2^{2k_{\xi}+1}f_{\xi})\Big(\frac{\eps}{2},\frac{\eps}{2}\Big)=\binom{2k_{\xi}}{k_{\xi}}\Big|\tilde{A}_{k_{\xi}}\Big(T-\frac{\eps}{2}\Big)^T\xi\Big|^2>0.
\end{equation}
We aim at proving that for all $\xi \in \mathbb{S}^{n-1}$ in the unit sphere, there exist positive constants $c_{\xi}>0$, $0< \tilde{\eps}_{\xi} \leq \frac{\eps}{2} $ and an open neighborhood $V_{\xi}$ of $\xi$ in $\mathbb{S}^{n-1}$ such that 
\begin{equation}\label{re4}
\forall \frac{\eps}{2} \leq t < \tau \leq \frac{\eps}{2}+\tilde{\eps}_{\xi}, \forall \eta \in V_{\xi}, \quad \int_{t}^{\tau}|A(T-s)^TR(t,s)^T\eta|^2ds \geq c_{\xi}(\tau-t)^{2k_{\xi}+1},
\end{equation}
By analogy with~\cite[Proposition~3.2]{joh}, we proceed by contradiction. If the assertion (\ref{re4}) does not hold, there exist sequences of positive real numbers $(t_l)_{l \geq 0}$, $(\tau_l)_{l \geq 0}$ satisfying
\begin{equation}\label{so0}
\forall l \geq 0, \quad \frac{\eps}{2} \leq t_l < \tau_l \leq \eps, \qquad \lim_{l \to +\infty}t_l=\lim_{l \to +\infty}\tau_l=\frac{\eps}{2},
\end{equation}
and a sequence $(\eta_l)_{l \geq 0}$ of elements in $\mathbb{S}^{n-1}$ so that 
\begin{equation}\label{re5}
\lim_{l \to +\infty}\eta_l=\xi, 
\end{equation}
and 
\begin{equation}\label{re5b}
\lim_{l \to +\infty}\frac{1}{(\tau_l-t_l)^{2k_{\xi}+1}}\int_{t_l}^{\tau_l}|A(T-s)^TR(t_l,s)^T\eta_l|^2ds=0.
\end{equation}
We deduce from (\ref{re5b}) that 
\begin{equation}\label{re6}
\lim_{l \to +\infty}\frac{1}{(\tau_l-t_l)^{2k_{\xi}+1}}\sup_{0 \leq t \leq \tau_l-t_l}\int_{t_l}^{t_l+t}|A(T-s)^TR(t_l,s)^T\eta_l|^2ds=0.
\end{equation}
Setting
\begin{equation}\label{re7}
u_l(x)=\frac{1}{(\tau_l-t_l)^{2k_{\xi}+1}}\int_{t_l}^{t_l+x(\tau_l-t_l)}|A(T-s)^TR(t_l,s)^T\eta_l|^2ds \geq 0, \quad 0 \leq x \leq 1,
\end{equation}
we can reformulate (\ref{re6}) as 
\begin{equation}\label{re8}
\lim_{l \to +\infty}\sup_{0 \leq x \leq 1}|u_l(x)|=0.
\end{equation}
By writing that
\begin{multline}\label{so1}
f_{\eta_l}(t_l,\tau)=\int_{t_l}^{\tau}|A(T-s)^TR(t_l,s)^T\eta_l|^2ds=\sum_{j=0}^{2k_{\xi}+1}a_l^{(j)}(\tau-t_l)^{j}\\
+\frac{(\tau-t_l)^{2k_{\xi}+2}}{(2k_{\xi}+1)!}\int_0^1(1-s)^{2k_{\xi}+1}(\partial_2^{2k_{\xi}+2}f_{\eta_l})\big(t_l,t_l+s(\tau-t_l)\big)ds,
\end{multline}
with $a_l^{(j)}=(\partial_2^jf_{\eta_l})(t_l,t_l)(j!)^{-1}$, we notice from (\ref{re1}) that there exists a positive constant $M>0$ such that 
\begin{multline}\label{so2}
\forall l \geq 0,\forall \tau \in [0,T], \\ \Big|\frac{1}{(2k_{\xi}+1)!}\int_0^1(1-s)^{2k_{\xi}+1}(\partial_2^{2k_{\xi}+2}f_{\eta_l})\big(t_l,t_l+s(\tau-t_l)\big)ds\Big| \leq M.
\end{multline}
It follows from (\ref{re7}), (\ref{so1}) and (\ref{so2}) that 
\begin{equation}\label{re9}
\forall 0 \leq x \leq 1, \forall l \geq 0, \quad \Big|u_l(x)-\sum_{j=0}^{2k_{\xi}+1}\frac{a_l^{(j)}}{(\tau_l-t_l)^{2k_{\xi}+1-j}}x^{j}\Big| \leq M(\tau_l-t_l) x^{2k_{\xi}+2}.
\end{equation}
It follows from (\ref{so0}), (\ref{re8}) and (\ref{re9}) that 
\begin{equation}\label{re10}
\lim_{l \to +\infty}\sup_{0 \leq x \leq 1}|p_l(x)|=0,
\end{equation}
with 
\begin{equation}\label{re11}
p_l(x)=\sum_{j=0}^{2k_{\xi}+1}\frac{a_l^{(j)}}{(\tau_l-t_l)^{2k_{\xi}+1-j}}x^{j}.
\end{equation}
By using the equivalence of norms in finite-dimensional vector space, we deduce from (\ref{re10}) that 
\begin{equation}\label{re12}
\forall 0 \leq j \leq 2k_{\xi}+1, \quad \lim_{l \to +\infty}\frac{a_l^{(j)}}{(\tau_l-t_l)^{2k_{\xi}+1-j}}=0.
\end{equation}
We obtain in particular that 
\begin{equation}\label{re13}
\lim_{l \to +\infty}a_l^{(2k_{\xi}+1)}=0.
\end{equation}
According to (\ref{re0b}), this is in contradiction with the fact that 
\begin{equation}\label{re14}
\lim_{l \to +\infty}a_l^{(2k_{\xi}+1)}=\lim_{l \to +\infty}\frac{(\partial_2^{2k_{\xi}+1}f_{\eta_l})(t_l,t_l)}{(2k_{\xi}+1)!}=\frac{1}{(2k_{\xi}+1)!}(\partial_2^{2k_{\xi}+1}f_{\xi})\Big(\frac{\eps}{2},\frac{\eps}{2}\Big)>0.
\end{equation}
By covering the compact set $\mathbb{S}^{n-1}$ by finitely many open neighborhoods of the form $(V_{\xi_j})_{1 \leq j \leq N}$, and letting $c=\inf_{1 \leq j \leq N}c_{\xi_j}>0$, $0<\tilde{\eps}=\inf_{1 \leq j \leq N}\tilde{\eps}_{\xi_j} \leq \eps/2 <1$, we conclude that 
$$\forall \xi \in \rr^n, \forall \frac{\eps}{2} \leq t \leq \tau \leq \frac{\eps}{2}+\tilde{\eps}, \quad \int_{t}^{\tau}|A(T-s)^TR(t,s)^T\xi|^2ds \geq c  (\tau-t)^{2n-1}|\xi|^2,$$
since $0 \leq k_{\xi} \leq n-1$.
This ends the proof of Lemma~\ref{alg}.
\end{proof}

In the following, we denote $\pi_j:L^2(\mathbb{R}^n) \rightarrow E_j$ the orthogonal frequency cutoff projection onto the closed subspace
\begin{equation}\label{sdf6}
E_j=\big\{ f \in L^2(\mathbb{R}^n) : \text{supp}(\hat{f}) \subset \{\xi \in \rr^n : |\xi| \leq j\}\big\}, \quad j \geq 0.
\end{equation}
Lemma~\ref{alg} allows to get the following exponential decay estimates of high frequencies of the solution to the Cauchy problem (\ref{adj_general1}):

\medskip

\begin{proposition} \label{thm:GS.z}
When the Kalman type condition (\ref{Kalman_time}) holds for some positive time $T>0$ belonging to the interval $I$, then there exists a positive constant $C>0$ such that 
for all $g_0 \in L^2(\rr^n)$, the mild solution $g$ to the Cauchy problem (\ref{adj_general1}) associated to the initial datum $g_0  \in L^2(\rr^n)$ satisfies
\begin{equation}\label{eq6.z}
\forall \frac{\eps}{2} \leq t \leq \tau \leq \frac{\eps}{2}+\tilde{\eps}, \forall k \geq 0, \quad  \| (1-\pi_k)g(\tau)\|_{L^2(\rr^n)} \leq e^{- C(\tau-t)^{2n-1}k^2} \|g(t)\|_{L^2(\rr^n)},
\end{equation}
with $0<\eps<T$ and $0<\tilde{\eps} \leq \eps/2$ being the positive constants defined in Lemma~\ref{alg}.
In particular, we have that for all $\frac{\eps}{2} \leq t_1<t_2<t_3<t_4 \leq \frac{\eps}{2}+\tilde{\eps} $, $k \geq 0$, $g_0 \in L^2(\rr^n)$, the mild solution $g$ to the Cauchy problem (\ref{adj_general1}) associated to the initial datum $g_0  \in L^2(\rr^n)$ satisfies
\begin{equation}\label{eq6_bis_1.z}
\|(1-\pi_k)g(t_3)\|_{L^2(\rr^n)}^2  \leq    e^{-C(t_3-t_2)^{2n-1}k^2 } \fint_{t_1}^{t_2} \|g(t)\|_{L^2(\rr^n)}^2 dt\,, 
\end{equation}
and
\begin{equation}\label{eq6_bis_2.z}
\fint_{t_3}^{t_4}\| (1-\pi_k)g(t)\|_{L^2(\rr^n)}^2 dt 
\leq    e^{- C(t_3-t_2)^{2n-1}k^2 } \fint_{t_1}^{t_2} \|g(t)\|_{L^2(\rr^n)}^2 dt\,,
\end{equation}
with the notation
$\fint_a^bf(t)dt=\frac{1}{b-a} \int_a^bf(t)dt.$
\end{proposition}

\medskip

\begin{proof}
Let $g_0 \in L^2(\mathbb{R}^n)$ and $g(t)$ be the mild solution of
$$\left\lbrace \begin{array}{l}
\partial_t g(t,x) - \frac{1}{2}\textrm{Tr}\big(A(T-t)A(T-t)^T\nabla_x^2 g(t,x)\big) + \big\langle B(T-t)x , \nabla_x g(t,x) \big\rangle\\
\hspace{9cm}  +\frac{1}{2}\textrm{Tr}\big(B(T-t)\big)g(t,x)= 0\,,\\
g|_{t=0}=g_0 \in L^2(\rr^n).
\end{array}\right.$$
We refer the reader to the appendix in Section~\ref{appendix} for the well-posedness of this Cauchy problem.
We deduce from the formula (\ref{sdf1}) in appendix that for all $0 \leq t \leq \tau \leq T$, 
\begin{multline}\label{sdf1.h}
\widehat{g}(\tau,\xi) 
= \exp\Big(\frac{1}{2}\int_{t}^{\tau}\textrm{Tr}\big(B(T-s)\big)ds\Big)\\
\times \widehat{g}\big(t,R(\tau,t)^T\xi\big)\exp\Big(-\frac{1}{2}\int_{t}^{\tau}|A(T-s)^TR(\tau,s)^T\xi|^2 ds\Big),
\end{multline}
where the resolvent $R$ is defined in (\ref{resol}).
It follows from (\ref{sdf1.h}) that for all $0 \leq t \leq \tau \leq T$, 
\begin{align}\label{sdf3.h}
& \  \|g(\tau)\|_{L^2(\rr^n)}^2\\ \notag
= & \ \frac{1}{(2\pi)^n}e^{\int_{t}^{\tau}\textrm{Tr}(B(T-s))ds}
\int_{\rr^n}\big|\widehat{g}\big(t,R(\tau,t)^T\xi\big)\big|^2e^{-\int_{t}^{\tau} |A(T-s)^TR(\tau,s)^T \xi|^2 ds}d\xi\\ \notag
= & \ \frac{\big|\text{det}\big(R(t,\tau)\big)\big|}{(2\pi)^n}e^{\int_{t}^{\tau}\textrm{Tr}(B(T-s))ds}
\int_{\rr^n}|\widehat{g}(t,\xi)|^2e^{-\int_{t}^{\tau} |A(T-s)^TR(t,s)^T \xi|^2 ds}d\xi\\ \notag
\leq & \ \frac{1}{(2\pi)^n}
\int_{\rr^n}|\widehat{g}(t,\xi)|^2d\xi= \|g(t)\|_{L^2(\rr^n)}^2,
\end{align}
since
\begin{equation}\label{liouville}
\forall t_1,t_2,t_3 \in [0,T], \quad R(t_1,t_2)R(t_2,t_3)=R(t_1,t_3)
\end{equation}
see e.g.~\cite{coron_book} (Proposition~1.5), and by Liouville formula 
\begin{equation}\label{liouville1}
\forall t, \tau \in [0,T], \quad \text{det}\big(R(\tau,t)\big)=\exp\Big(\int_{t}^{\tau}\textrm{Tr}\big(B(T-s)\big)ds\Big),
\end{equation}
see e.g.~\cite{cartan} (Proposition II.2.3.1).
This implies that the function $t \mapsto \|g(t)\|_{L^2(\rr^n)}$ is non-increasing
\begin{equation}\label{estt200.z}
\forall 0 \leq t \leq \tau \leq T, \quad \|g(\tau)\|_{L^2(\rr^n)} \leq \|g(t)\|_{L^2(\rr^n)}.
\end{equation}
We deduce from (\ref{sdf1.h}), (\ref{liouville}), (\ref{liouville1}) and Lemma~\ref{alg} that for all $\frac{\eps}{2} \leq t \leq \tau \leq \frac{\eps}{2}+\tilde{\eps}$, $k \geq 0$, $g_0 \in L^2(\rr^n)$, the mild solution $g$ to the Cauchy problem (\ref{adj_general1}) associated to the initial datum $g_0  \in L^2(\rr^n)$ satisfies
\begin{align}\label{sdf3}
& \  \| (1-\pi_k)g(\tau)\|_{L^2(\rr^n)}^2\\ \notag
= & \ \frac{1}{(2\pi)^n}e^{\int_t^{\tau}\textrm{Tr}(B(T-s))ds}
\int_{|\xi| \geq k}\big|\widehat{g}\big(t,R(\tau,t)^T\xi\big)\big|^2e^{-\int_t^{\tau} |A(T-s)^TR(\tau,s)^T \xi|^2 ds}d\xi\\ \notag
= & \ \frac{1}{(2\pi)^n}\big|\text{det}\big(R(t,\tau)\big)\big|e^{\int_t^{\tau}\textrm{Tr}(B(T-s))ds}
\int_{|R(t,\tau)^T\xi| \geq k}|\widehat{g}(t,\xi)|^2e^{-\int_t^{\tau} |A(T-s)^TR(t,s)^T \xi|^2 ds}d\xi\\ \notag
\leq & \ \frac{1}{(2\pi)^n}
\int_{|\xi| \geq k |R(t,\tau)^T|^{-1}}|\widehat{g}(t,\xi)|^2e^{-c(\tau-t)^{2n-1}|\xi|^2}d\xi.
\end{align}
We deduce from (\ref{sdf3}) that there exists $C>0$ such that for all $\frac{\eps}{2} \leq t \leq \tau \leq \frac{\eps}{2}+\tilde{\eps}$, $k \geq 0$, $g_0 \in L^2(\rr^n)$, the mild solution $g$ to the Cauchy problem (\ref{adj_general1}) associated to the initial datum $g_0  \in L^2(\rr^n)$ satisfies
\begin{equation}\label{sdf5}
\| (1-\pi_k)g(\tau)\|_{L^2(\rr^n)}^2 \leq e^{-2C(\tau-t)^{2n-1}k^2}\|g(t)\|_{L^2(\rr^n)}^2.
\end{equation}
It proves the estimate (\ref{eq6.z}). We therefore deduce (\ref{eq6_bis_1.z}) from (\ref{eq6.z}) and (\ref{estt200.z}), that is, that for all $\frac{\eps}{2} \leq t_1 < t_2 <t_3 \leq \frac{\eps}{2}+\tilde{\eps}$, $k \geq 0$, $g_0 \in L^2(\rr^n)$, the mild solution $g$ to the Cauchy problem (\ref{adj_general1}) associated to the initial datum $g_0  \in L^2(\rr^n)$ satisfies
\begin{multline*}
\|(1-\pi_k)g(t_3)\|_{L^2(\rr^n)}^2  \leq e^{-2C(t_3-t_2)^{2n-1}k^2}\|g(t_2)\|_{L^2(\rr^n)}^2 \\
 \leq  e^{-2C(t_3-t_2)^{2n-1}k^2} \fint_{t_1}^{t_2} \|g(t)\|_{L^2(\rr^n)}^2 dt.
 \end{multline*}
By using anew (\ref{eq6.z}), it follows that for all $\frac{\eps}{2} \leq t_1<t_2<t_3<t_4 \leq \frac{\eps}{2}+\tilde{\eps}$, $k \geq 0$, $g_0 \in L^2(\rr^n)$,
the mild solution $g$ to the Cauchy problem (\ref{adj_general1}) associated to the initial datum $g_0  \in L^2(\rr^n)$ satisfies
\begin{align*}
& \ \fint_{t_3}^{t_4}\| (1-\pi_k)g(t)\|_{L^2(\rr^n)}^2dt \\
\leq &\ \fint_{t_3}^{t_4}\exp\Big(-2C\Big(t-t_1-\frac{t-t_3}{t_4-t_3}(t_2-t_1)\Big)^{2n-1}k^2\Big)\Big\|g\Big(t_1+\frac{t-t_3}{t_4-t_3}(t_2-t_1)\Big)\Big\|_{L^2(\rr^n)}^2dt \\
\leq & \ e^{-2C(t_3-t_2)^{2n-1}k^2}\fint_{t_3}^{t_4}\Big\|g\Big(t_1+\frac{t-t_3}{t_4-t_3}(t_2-t_1)\Big)\Big\|_{L^2(\rr^n)}^2dt\\
= & \ e^{-2C(t_3-t_2)^{2n-1}k^2}\fint_{t_1}^{t_2}\|g(t)\|_{L^2(\rr^n)}^2dt, 
\end{align*}
since 
$$\forall t_3 \leq t \leq t_4, \quad t-t_1-\frac{t-t_3}{t_4-t_3}(t_2-t_1) \geq t_3-t_2>0.$$
This ends the proof of Proposition~\ref{thm:GS.z}.
\end{proof}

\medskip

Before pursuing the proof of Theorem~\ref{meta_thm}, we notice by using the very same lines as in the proof of Proposition~\ref{thm:GS.z} that the following Gevrey $1/2$ type estimates hold:

\medskip

\begin{corollary}\label{corollary}
Under the assumptions of Proposition~\ref{thm:GS.z}, and with $0<\eps<T$ and $0<\tilde{\eps} \leq \eps/2$ the positive constants defined in Proposition~\ref{thm:GS.z}, then there exists a positive constant $c_0>0$ such that for all $g_0 \in L^2(\rr^n)$, the mild solution $g$ to the Cauchy problem (\ref{adj_general1}) associated to the initial datum $g_0  \in L^2(\rr^n)$ satisfies
\begin{multline}\label{ser1}
\forall \frac{\eps}{2} \leq t < \tau \leq \frac{\eps}{2}+\tilde{\eps}, \forall k \geq 0, \\  \||D_x|^{k}g(\tau)\|_{L^2(\rr^n)}  \leq \frac{c_0^k}{(\tau-t)^{\frac{(2n-1)k}{2}}}\sqrt{k!}\|g(t)\|_{L^2(\rr^n)}\leq \frac{c_0^k}{(\tau-t)^{\frac{(2n-1)k}{2}}}\sqrt{k!}\|g_0\|_{L^2(\rr^n)}.
\end{multline}
\end{corollary}

\medskip

\begin{proof}
We deduce from (\ref{sdf1.h}), (\ref{liouville}), (\ref{liouville1}) and Lemma~\ref{alg} that for all $\frac{\eps}{2} \leq t < \tau \leq \frac{\eps}{2}+\tilde{\eps}$, $k \geq 0$, $g_0 \in L^2(\rr^n)$, the mild solution $g$ to the Cauchy problem (\ref{adj_general1}) associated to the initial datum $g_0  \in L^2(\rr^n)$ satisfies
\begin{align}\label{sdf3.m}
& \  \||D_x|^{k}g(\tau)\|_{L^2(\rr^n)}^2\\ \notag
= & \ \frac{1}{(2\pi)^n}e^{\int_t^{\tau}\textrm{Tr}(B(T-s))ds}
\int_{\rr^n}|\xi|^{2k}\big|\widehat{g}\big(t,R(\tau,t)^T\xi\big)\big|^2e^{-\int_t^{\tau} |A(T-s)^TR(\tau,s)^T \xi|^2 ds}d\xi\\ \notag
= & \ \frac{1}{(2\pi)^n}\big|\text{det}\big(R(t,\tau)\big)\big|e^{\int_t^{\tau}\textrm{Tr}(B(T-s))ds}
\int_{\rr^n}|R(t,\tau)^T\xi|^{2k}|\widehat{g}(t,\xi)|^2e^{-\int_t^{\tau} |A(T-s)^TR(t,s)^T \xi|^2 ds}d\xi\\ \notag
\leq & \ \frac{|R(t,\tau)|^{2k}}{(2\pi)^n}
\int_{\rr^n}|\xi|^{2k}|\widehat{g}(t,\xi)|^2e^{-c(\tau-t)^{2n-1}|\xi|^2}d\xi.
\end{align}
According to (\ref{estt200.z}), this implies that 
\begin{align*}
 & \ \||D_x|^{k}g(\tau)\|_{L^2(\rr^n)}^2 \\
 \leq & \  \frac{|R(t,\tau)|^{2k}k!}{(2\pi)^nc^k(\tau-t)^{(2n-1)k}}
\int_{\rr^n}\frac{c^k(\tau-t)^{(2n-1)k}|\xi|^{2k}}{k!}|\widehat{g}(t,\xi)|^2e^{-c(\tau-t)^{2n-1}|\xi|^2}d\xi\\
\leq & \  \frac{|R(t,\tau)|^{2k}k!}{c^k(\tau-t)^{(2n-1)k}}\|g(t)\|_{L^2(\rr^n)}^2
\leq  \frac{|R(t,\tau)|^{2k}k!}{c^k(\tau-t)^{(2n-1)k}}\|g_0\|_{L^2(\rr^n)}^2.
\end{align*}
This ends the proof of Corollary~\ref{corollary}.
\end{proof}

\medskip

We resume the proof of Theorem~\ref{meta_thm} by adapting the Lebeau-Robbiano direct approach for observability. To that end, we use the results of exponential decay given by Proposition~\ref{thm:GS.z} and the following spectral inequality proved by Le Rousseau and Moyano in \cite[Theorem 3.1]{LeRousseau_Moyano}:

\medskip

\begin{theorem} \label{thm:IS_LR_M}
If $\omega$ is an open subset of $\mathbb{R}^n$ satisfying condition (\ref{hyp_omega}), then there exists a positive constant $c_1>1$ such that 
$$\|f\|_{L^2(\mathbb{R}^n)} \leqslant c_1e^{c_1 N} \|f\|_{L^2(\omega)},$$
for all $N \geq 0$ and $f \in L^2(\mathbb{R}^n)$ whose Fourier transform verifies 
$$\emph{\text{supp}}(\hat{f})\subset \{\xi \in \rr^n : |\xi| \leq N\},$$ 
with $|\cdot|$ being the Euclidean norm on $\rr^n$. 
\end{theorem}

\medskip

\noindent
We can now provide the core of the proof of Theorem~\ref{meta_thm}:

\medskip

\noindent
\textit{Step 1: Preliminaries.} Let  $\rho$ be a positive constant satisfying 
\begin{equation}\label{eq8.z}
0 < \rho <  \frac{1}{2n-1}.
\end{equation} 
We consider $K>0$ the positive constant verifying 
$$\sum_{k=1}^{+\infty}\frac{2K}{4^{k\rho}}=\tilde{\eps},$$
where $\tilde{\eps}>0$ is the positive constant defined in Lemma~\ref{alg}. 
We define for all $k \geq 1$,
\begin{equation}\label{eq8.1.z}
\tau_k=\frac{K}{4^{k \rho}}, \quad \alpha_0=0, \quad \alpha_k=\sum_{j=1}^k2\tau_j, \quad J_k=\Big[\frac{\eps}{2}+\tilde{\eps}-\alpha_{k-1}-\tau_k,\frac{\eps}{2}+\tilde{\eps}-\alpha_{k-1}\Big],
\end{equation}
where $\eps>0$ is the positive constant defined in Lemma~\ref{alg}.
We observe that the sequence $(\alpha_k)_{k \geq 0}$ is increasing and that $\lim_{k \to +\infty}\alpha_k=\tilde{\eps}$.
\begin{figure}[h]
\begin{picture}(450,40)
\put(218,-10){\line(0,1){35}} 
\put(170,25){\line(1,0){48}} 
\put(170,25){\line(0,-1){20}} 
\put(92,-5){$\ \ \frac{\eps}{2}+\tilde{\eps}-\alpha_k$}
\put(220,-5){$\frac{\eps}{2}+\tilde{\eps}-\alpha_{k-1}$} 
\put(185,30){$J_k$}
\put(170,-10){$2\tau_k$}
\put(125,-10){\line(0,1){25}}
\thicklines 
\put(120,5){\line(1,0){110}}
\end{picture}
\end{figure}
According to (\ref{eq8.z}), we can choose $\beta>0$ a positive constant satisfying
\begin{equation}\label{def:alpha.z}
1+\rho(2n-1) < \beta < 2\,.
\end{equation}
We define for all $k \geq 1$,
\begin{equation} \label{def:lk.z}
l_k=[2^{k \beta}],
\end{equation}
the integer part of $2^{k \beta}$.
We claim that there exists an integer $p_0 \geq 2$ such that
\begin{equation} \label{def_p0_1.z}
\forall k \geq p_0\,, \quad  e^{ - C(4^{k-1}\alpha_{k-1}^{2n-1} + l_k^2 \tau_k^{2n-1} - 4^k \alpha_k^{2n-1})} \leqslant \frac{1}{2}\,,
\end{equation}
\begin{equation} \label{def_p0_2.z}
\forall k \geq p_0\,, \quad 2 c_1 l_k - C4^{k-1}\alpha_{k-1}^{2n-1} \leqslant 0\,,
\end{equation}
\begin{equation} \label{def_p0_4.z}
\forall k \geq p_0\,, \quad \frac{e^{ 2 c_1 l_k - C4^{k-2}\alpha_{k-1}^{2n-1}}}{\tau_k} \leqslant 1,
\end{equation}
\begin{equation} \label{def_p0_3.z}
\forall k \geq p_0\,, \quad 2 c_1^2 e^{ 2 c_1 l_k - C l_k^2 \tau_k^{2n-1}} \leqslant 1,
\end{equation}
where the positive constant $C>0$ is defined in Proposition~\ref{thm:GS.z}, whereas the positive constant $c_1>1$ is defined in Theorem~\ref{thm:IS_LR_M}.
Indeed, the claim (\ref{def_p0_1.z}) follows from (\ref{eq8.1.z}), (\ref{def:alpha.z}) and (\ref{def:lk.z}) as
\begin{equation}\label{estt1.z}
l_k^2 \tau_k^{2n-1} \underset{k \rightarrow +\infty}{\sim} K^{2n-1} 4^{k[\beta-\rho(2n-1)]}\,,  \ 
 4^{k-1}\alpha_{k-1}^{2n-1} - 4^k \alpha_k^{2n-1} \underset{k \rightarrow +\infty}{\sim} -4^{k-1}(3\tilde{\eps}^{2n-1}).
\end{equation}
and $\beta-\rho(2n-1) > 1$. The claims (\ref{def_p0_2.z}) and (\ref{def_p0_4.z}) follow from (\ref{eq8.1.z}), (\ref{def:alpha.z}) and (\ref{def:lk.z}) as
\begin{equation}\label{estt2.z}
4^{k-1} \alpha_{k-1}^{2n-1} \underset{k \rightarrow +\infty}{\sim} 4^{k-1} \tilde{\eps}^{2n-1}\,, \qquad
l_k \underset{k \rightarrow +\infty}{\sim} 2^{\beta k}, \qquad \beta<2.
\end{equation}
Finally, the claim (\ref{def_p0_3.z}) follows from (\ref{eq8.z}), (\ref{def:alpha.z}), (\ref{estt1.z}) and (\ref{estt2.z}) as
$$2\beta-2\rho(2n-1) > \beta.$$

\bigskip

\noindent
\textit{Step 2.} In this second step, we aim at establishing that for all $g_0 \in L^2(\rr^n)$,
the mild solution $g$ to the Cauchy problem (\ref{adj_general1}) associated to the initial datum $g_0  \in L^2(\rr^n)$ satisfies the following estimate:
\begin{multline}\label{eq9.4.z}
\|g(T)\|_{L^2(\rr^n)}^2 \leq \Big\|g\Big(\frac{\eps}{2}+\tilde{\eps}\Big)\Big\|_{L^2(\rr^n)}^2\\
 \leq \Big\|\pi_{2^{p_0-1}}g\Big(\frac{\eps}{2}+\tilde{\eps}\Big)\Big\|_{L^2(\rr^n)}^2+ 2 \sum_{k=p_0}^{+\infty} e^{-C4^{k-1}\alpha_{k-1}^{2n-1}} \fint_{J_k}\|\pi_{l_k}g(t)\|_{L^2(\rr^n)}^2dt.
\end{multline}
To that end, we notice that the first inequality follows from Lemma~\ref{alg} and (\ref{estt200.z}) since $\frac{\eps}{2}+\tilde{\eps} \leq \eps \leq T$, and we deduce from (\ref{sdf6}) and (\ref{eq6_bis_1.z}) that
\begin{align}\label{eq9.00.z}
\|g(T)\|_{L^2(\rr^n)}^2 \leq & \  \Big\|g\Big(\frac{\eps}{2}+\tilde{\eps}\Big)\Big\|_{L^2(\rr^n)}^2 \notag \\ 
= &\ \Big\|\pi_{2^{p_0-1}}g\Big(\frac{\eps}{2}+\tilde{\eps}\Big)\Big\|_{L^2(\rr^n)}^2 + \sum_{k=p_0}^{+\infty}\Big\|(\pi_{2^k}-\pi_{2^{k-1}})g\Big(\frac{\eps}{2}+\tilde{\eps}\Big)\Big\|_{L^2(\rr^n)}^2\notag \\  \notag
\leq &\  \Big\|\pi_{2^{p_0-1}}g\Big(\frac{\eps}{2}+\tilde{\eps}\Big)\Big\|_{L^2(\rr^n)}^2 + \sum_{k=p_0}^{+\infty}\Big\| (1-\pi_{2^{k-1}})g\Big(\frac{\eps}{2}+\tilde{\eps}\Big)\Big\|_{L^2(\rr^n)}^2 \\ 
\leq &\ \Big\|\pi_{2^{p_0-1}}g\Big(\frac{\eps}{2}+\tilde{\eps}\Big)\Big\|_{L^2(\rr^n)}^2+ \sum_{k=p_0}^{+\infty} e^{-C4^{k-1}\alpha_{k-1}^{2n-1}} \fint_{J_k}\|g(t)\|_{L^2(\rr^n)}^2 dt.
\end{align}
By writing that 
\begin{multline*}
\sum_{k=p_0}^{+\infty} e^{-C4^{k-1}\alpha_{k-1}^{2n-1}} \fint_{J_k}\|g(t)\|_{L^2(\rr^n)}^2 dt \\
= \sum_{k=p_0}^{+\infty} e^{-C4^{k-1}\alpha_{k-1}^{2n-1}} \fint_{J_k}(\|\pi_{l_k}g(t)\|_{L^2(\rr^n)}^2+\|(1-\pi_{l_k})g(t)\|_{L^2(\rr^n)}^2) dt,
\end{multline*}
it follows from (\ref{eq6_bis_2.z}) and (\ref{def_p0_1.z}) that
\begin{align*}
& \sum_{k=p_0}^{+\infty} e^{-C4^{k-1}\alpha_{k-1}^{2n-1}} \fint_{J_k}\|g(t)\|_{L^2(\rr^n)}^2 dt \leq \ \sum_{k=p_0}^{+\infty} e^{-C4^{k-1}\alpha_{k-1}^{2n-1}} \fint_{J_k}  \|\pi_{l_k}g(t)\|_{L^2(\rr^n)}^2 dt\\ \notag
& \qquad + \sum_{k=p_0}^{+\infty} e^{-C(4^{k-1}\alpha_{k-1}^{2n-1} + l_k^2 \tau_k^{2n-1}) } \fint_{J_{k+1}}\|g(t)\|_{L^2(\rr^n)}^2 dt \\ \notag
\leq &\ \sum_{k=p_0}^{+\infty} e^{-C4^{k-1}\alpha_{k-1}^{2n-1}} \fint_{J_k}\|\pi_{l_k}g(t)\|_{L^2(\rr^n)}^2 dt
+ \frac{1}{2} \sum_{k=p_0}^{+\infty} e^{-C4^{k}\alpha_{k}^{2n-1}} \fint_{J_{k+1}}\|g(t)\|_{L^2(\rr^n)}^2 dt,
\end{align*}
implying that 
$$\sum_{k=p_0}^{+\infty} e^{-C4^{k-1}\alpha_{k-1}^{2n-1}} \fint_{J_k}\|g(t)\|_{L^2(\rr^n)}^2 dt
\leqslant 2 \sum_{k=p_0}^{+\infty} e^{-C4^{k-1}\alpha_{k-1}^{2n-1}} \fint_{J_k}\|\pi_{l_k}g(t)\|_{L^2(\rr^n)}^2 dt.$$
The previous estimate together with (\ref{eq9.00.z}) provides (\ref{eq9.4.z}).

\bigskip

\noindent
\textit{Step 3: Induction.} We prove by induction that for all $N \geq p_0$, $g_0 \in L^2(\rr^n)$, the mild solution $g$ to the Cauchy problem (\ref{adj_general1}) associated to the initial datum $g_0  \in L^2(\rr^n)$ satisfies
\begin{multline}\label{eq9.5.z}
\Big\|\pi_{2^{p_0-1}}g\Big(\frac{\eps}{2}+\tilde{\eps}\Big)\Big\|_{L^2(\rr^n)}^2+2  \sum_{k=p_0}^{N} e^{-C4^{k-1}\alpha_{k-1}^{2n-1}} \fint_{J_k}\|\pi_{l_k}g(t)\|_{L^2(\rr^n)}^2dt\\
\leq \sum_{k=p_0}^{N} a_k  \fint_{J_k}\|g(t)\|_{L^2(\omega)}^2dt+B_N \fint_{J_{N+1}}\|g(t)\|_{L^2(\rr^n)}^2 dt,
\end{multline}
with
\begin{equation} \label{def:aN.z}
a_{p_0}= 4 c_1^2, \quad B_{p_0}= 1 +  4  c_1^2,  
\end{equation}
\begin{equation} \label{def:BN.z}
\forall N \geq p_0, \quad a_{N+1}= 2 c_1^2 e^{2 c_1 l_{N+1}}(2e^{-C4^{N}\alpha_{N}^{2n-1}}+B_N),
\end{equation}\begin{equation} \label{def:BN1.z}
\forall N \geq p_0, \quad B_{N+1}= \left( a_{N+1}+B_N \right) e^{- Cl_{N+1}^2 \tau_{N+1}^{2n-1}}.
\end{equation}
\textit{Initialization:} We observe from (\ref{estt200.z}) and (\ref{eq8.1.z}) that
\begin{equation} \label{init_1.z}
\Big\|\pi_{2^{p_0-1}}g\Big(\frac{\eps}{2}+\tilde{\eps}\Big)\Big\|_{L^2(\rr^n)}^2 \leq \Big\|g\Big(\frac{\eps}{2}+\tilde{\eps}\Big)\Big\|_{L^2(\rr^n)}^2 \leq  \fint_{J_{p_0+1}}\|g(t)\|_{L^2(\rr^n)}^2 dt.
\end{equation}
We deduce from the spectral inequality (Theorem~\ref{thm:IS_LR_M}), (\ref{eq6_bis_2.z}) in Proposition~\ref{thm:GS.z} and (\ref{eq8.1.z}) that
\begin{align*}
& \quad  \fint_{J_{p_0}}\|\pi_{l_{p_0}}g(t)\|_{L^2(\rr^n)}^2dt \leq c_1^2 e^{2c_1 l_{p_0}} \fint_{J_{p_0}}\|\pi_{l_{p_0}}g(t)\|_{L^2(\omega)}^2dt \\ 
&\ \leq 2 c_1^2 e^{2 c_1 l_{p_0}} \fint_{J_{p_0}}(\|g(t)\|_{L^2(\omega)}^2+\|(1-\pi_{l_{p_0}})g(t)\|_{L^2(\rr^n)}^2) dt \\  
&\ \leq 2 c_1^2 e^{2 c_1 l_{p_0}} \fint_{J_{p_0}}\|g(t)\|_{L^2(\omega)}^2dt + 2  c_1^2 e^{2 c_1 l_{p_0}-Cl_{p_0}^2\tau_{p_0}^{2n-1}} \fint_{J_{p_0+1}}\|g(t)\|_{L^2(\rr^n)}^2dt,
\end{align*}
since $\|(1-\pi_{l_{p_0}})g(t)\|_{L^2(\omega)} \leq \|(1-\pi_{l_{p_0}})g(t)\|_{L^2(\rr^n)}$.
By summing the two previous inequalities, we obtain from (\ref{def_p0_2.z}) that the initialization of the induction process holds:
\begin{align}
&\ \Big\|\pi_{2^{p_0-1}}g\Big(\frac{\eps}{2}+\tilde{\eps}\Big)\Big\|_{L^2(\rr^n)}^2+2e^{-C4^{p_0-1}\alpha_{p_0-1}^{2n-1}}\fint_{J_{p_0}}\|\pi_{l_{p_0}}g(t)\|_{L^2(\rr^n)}^2dt \\ \notag
\leq  &\ 4c_1^2e^{2c_1 l_{p_0}-C4^{p_0-1}\alpha_{p_0-1}^{2n-1}}\fint_{J_{p_0}}\|g(t)\|_{L^2(\omega)}^2 dt \\ \notag
&\ +(1+4 c_1^2e^{2c_1 l_{p_0}-C4^{p_0-1}\alpha_{p_0-1}^{2n-1}})\fint_{J_{p_0+1}}\|g(t)\|_{L^2(\rr^n)}^2dt\\ \notag
\leq  &\ 4c_1^2\fint_{J_{p_0}}\|g(t)\|_{L^2(\omega)}^2 dt +(1+4 c_1^2)\fint_{J_{p_0+1}}\|g(t)\|_{L^2(\rr^n)}^2dt.
\end{align}
\textit{Heredity:} Let us assume that (\ref{eq9.5.z}) holds at the rank $N \geq p_0$. 
We deduce from (\ref{eq6_bis_2.z}) in Proposition~\ref{thm:GS.z} and (\ref{eq8.1.z}) that
\begin{align}\label{eq13.z}
& \ B_N \fint_{J_{N+1}}\|g(t)\|_{L^2(\rr^n)}^2 dt \\ \notag
= & \ B_N \fint_{J_{N+1}}(\|\pi_{l_{N+1}}g(t)\|_{L^2(\rr^n)}^2+\|(1-\pi_{l_{N+1}})g(t)\|_{L^2(\rr^n)}^2)dt\\ \notag
\leq & \ B_N \fint_{J_{N+1}}\|\pi_{l_{N+1}}g(t)\|_{L^2(\rr^n)}^2dt+B_N e^{- C l_{N+1}^2 \tau_{N+1}^{2n-1}} \fint_{J_{N+2}}\|g(t)\|_{L^2(\rr^n)}^2dt.
\end{align}
By using successively the induction assumption, the previous inequality, the spectral inequality (Theorem~\ref{thm:IS_LR_M}), the definition (\ref{def:BN.z}), the estimate (\ref{eq6_bis_2.z}) and the definition (\ref{def:BN1.z}), we obtain that
\begin{align}
     &\ \Big\|\pi_{2^{p_0-1}}g\Big(\frac{\eps}{2}+\tilde{\eps}\Big)\Big\|_{L^2(\rr^n)}^2+2 \sum_{k=p_0}^{N+1} e^{-C4^{k-1}\alpha_{k-1}^{2n-1}} \fint_{J_k}\|\pi_{l_k}g(t)\|_{L^2(\rr^n)}^2dt \\ \notag
\leq &\ \sum_{k=p_0}^{N} a_k   \fint_{J_k}\|g(t)\|_{L^2(\omega)}^2dt \\ \notag
     &\ + c_1^2 e^{2 c_1 l_{N+1}}(2e^{-C 4^{N}\alpha_{N}^{2n-1}}+B_N)\fint_{J_{N+1}}\|\pi_{l_{N+1}}g(t)\|_{L^2(\omega)}^2dt  \\ \notag
     &\  +  B_N e^{- C l_{N+1}^2 \tau_{N+1}^{2n-1}} \fint_{J_{N+2}}\|g(t)\|_{L^2(\rr^n)}^2dt \\ \notag
\leq &\ \sum_{k=p_0}^{N+1} a_k\fint_{J_k}\|g(t)\|_{L^2(\omega)}^2dt \\ \notag
     &\ + a_{N+1}\fint_{J_{N+1}}\|(1-\pi_{l_{N+1}})g(t)\|_{L^2(\rr^n)}^2dt  \\ \notag
     &\ +  B_Ne^{- C l_{N+1}^2 \tau_{N+1}^{2n-1}} \fint_{J_{N+2}}\|g(t)\|_{L^2(\rr^n)}^2dt \\ \notag
\leq  &\ \sum_{k=p_0}^{N+1} a_k  \fint_{J_k}\|g(t)\|_{L^2(\omega)}^2dt+B_{N+1} \fint_{J_{N+2}}\|g(t)\|_{L^2(\rr^n)}^2dt,
\end{align}
since $\|(1-\pi_{l_{N+1}})g(t)\|_{L^2(\omega)} \leq \|(1-\pi_{l_{N+1}})g(t)\|_{L^2(\rr^n)}$.
This proves that the estimate (\ref{eq9.5.z}) holds for all $N \geq p_0$.

\bigskip

\noindent
\textit{Step 4.} The aim of this step is to prove that the sequence $(a_N \tau_N^{-1})_{N \geq p_0}$ is bounded and that $\lim_{N \to +\infty}B_N=0$.
By using (\ref{def_p0_3.z}), (\ref{def:BN.z}) and (\ref{def:BN1.z}), we obtain that for all $N \geq p_0$,
\begin{multline} \label{Majo:B(N+1).z}
0 \leq B_{N+1}=  B_N\big(e^{-C l_{N+1}^2\tau_{N+1}^{2n-1}} + 2 c_1^2 e^{2c_1 l_{N+1}-C l_{N+1}^2\tau_{N+1}^{2n-1}}\big) \\ 
+ 2 \big(2c_1^2 e^{2c_1 l_{N+1}-C l_{N+1}^2\tau_{N+1}^{2n-1}}\big) e^{-C 4^N \alpha_N^{2n-1}} \leq 2 B_N + 2 .
\end{multline}
It follows from (\ref{Majo:B(N+1).z}) that 
$$\forall N \geq p_0, \quad 0 \leq B_{N+1} + 2 \leqslant 2  \left( B_N + 2 \right),$$
implying that for all $N \geqslant p_0$,
\begin{equation}\label{estt6.z}
0 \leqslant B_N \leqslant B_N+2 \leqslant  \left( B_{p_0} + 2 \right)2^{N-p_0}.
\end{equation}
By implementing (\ref{estt6.z}) in the estimate (\ref{Majo:B(N+1).z}), we obtain that there exists a positive constant $c(p_0,c_1)>0$ such that for all $N \geq p_0$,
\begin{equation}\label{estt7.z}
0 \leqslant B_{N+1} \leqslant c(p_0,c_1)\big(2^{N} e^{2c_1 l_{N+1}-C l_{N+1}^2 \tau_{N+1}^{2n-1}} +  e^{2 c_1 l_{N+1}-C 4^N \alpha_N^{2n-1}}\big).
\end{equation}
On the other hand, we observe from (\ref{eq8.1.z}), (\ref{def:alpha.z}) and (\ref{def:lk.z}) that 
\begin{equation}\label{estt8.z}
l_{N+1}^2\tau_{N+1}^{2n-1} \underset{N \rightarrow +\infty}{\sim} K^{2n-1} 4^{(N+1)[\beta-\rho(2n-1)]}\,, \quad
l_{N+1} \underset{N \rightarrow +\infty}{\sim} 2^{(N+1)\beta}
\end{equation}
and
\begin{equation}\label{estt9.z}
\beta<2<2\beta-2\rho(2n-1), \quad 4^N \alpha_N^{2n-1} \underset{N \rightarrow +\infty}{\sim}  4^{N} \tilde{\eps}^{2n-1}.
\end{equation}
It follows from (\ref{estt7.z}), (\ref{estt8.z}) and (\ref{estt9.z}) that 
\begin{equation}\label{estt10.z}
0 \leqslant B_N \leqslant e^{- C4^{N-1} \alpha_N^{2n-1}},
\end{equation}
when $N \gg 1$.
It proves that $\lim_{N \to +\infty} B_N=0$.
We deduce from (\ref{def_p0_4.z}), (\ref{def:BN.z}) and (\ref{estt10.z}) that 
$$0 \leq a_{N+1}\tau_{N+1}^{-1} = 4 c_1^2\frac{e^{2 c_1 l_{N+1}-C 4^N \alpha_N^{2n-1}}}{\tau_{N+1}} + 2 c_1^2 \frac{e^{2 c_1 l_{N+1}} B_N}{\tau_{N+1}}  
\leq 6c_1^2\frac{e^{2 c_1 l_{N+1}-C 4^{N-1}\alpha_N^{2n-1}}}{\tau_{N+1}} \leq 6c_1^2,$$
when $N \gg 1$. It shows that the sequence $(a_N\tau_N^{-1})_{N \geqslant p_0}$ is bounded.

\bigskip

\noindent
\textit{Step 5: Conclusion.} Let $M>0$ be a positive constant satisfying
$$\forall N \geq p_0, \quad 0 \leqslant a_N\tau_N^{-1} \leqslant M.$$ 
We deduce from (\ref{eq9.4.z}), (\ref{eq9.5.z}) while passing to the limit $N \rightarrow +\infty$ thanks to Step 4 that
\begin{multline}\label{estt11.z}
\|g(T)\|_{L^2(\rr^n)}^2  \leq  \Big\|\pi_{2^{p_0-1}}g\Big(\frac{\eps}{2}+\tilde{\eps}\Big)\Big\|_{L^2(\rr^n)}^2\\
+2 \sum_{k=p_0}^{+\infty} e^{-C 4^{k-1}\alpha_{k-1}^{2n-1}} \fint_{J_k}\|\pi_{l_k}g(t)\|_{L^2(\rr^n)}^2 dt \leq  M \sum_{k=p_0}^{+\infty}  \int_{J_k}\|g(t)\|_{L^2(\omega)}^2 dt, 
\end{multline}
since
$$0 \leq \fint_{J_{N+1}}\|g(t)\|_{L^2(\rr^n)}^2 dt \leq \|g_0\|_{L^2(\rr^n)}^2,$$
according to (\ref{estt200.z}).
Since the intervals $(J_k)_{k \geqslant p_0}$ are disjoint and $\sqcup_{k=p_0}^{+\infty}J_k \subset [\frac{\eps}{2},\frac{\eps}{2}+\tilde{\eps}] \subset [0,T]$, we finally obtain that for all $g_0 \in L^2(\rr^n)$, the mild solution $g$ to the Cauchy problem (\ref{adj_general1}) associated to the initial datum $g_0  \in L^2(\rr^n)$ satisfies
$$\|g(T)\|_{L^2(\rr^n)}^2 \leqslant M \int_0^T\|g(t)\|_{L^2(\omega)}^2dt.$$
It proves the observability estimate (\ref{dfg2}) and ends the proof of Theorem \ref{meta_thm}\,.

\section{Appendix}\label{appendix}

\subsection{Well-posedness of the homogeneous and inhomogeneous Cauchy problems}
We first study the well-posedness of the homogeneous equation
\begin{equation} \label{CY_Hom}
\left\lbrace \begin{array}{l}
\partial_t k(t,x) - P_0(t) k(t,x)=0, \quad  (t,x) \in (t_0,T_1) \times \mathbb{R}^n, \\
k|_{t=t_0} = k_0 \in L^2(\mathbb{R}^n),
\end{array} \right.
\end{equation}
associated to the non-autonomous Ornstein-Uhlenbeck operator
\begin{equation} \label{def:P_appendix}
P_0(t)=\frac{1}{2} \text{Tr}\big(A_0(t)A_0(t)^T\nabla_x^2\big)-\big\langle B_0(t)x,\nabla_x\big\rangle-\frac{1}{2}\text{Tr}\big(B_0(t)\big), \quad t \in (T_0,T_1),
\end{equation}
with $T_0 \leqslant t_0 < T_1$ and $A_0, B_0 \in C^0([T_0,T_1],M_n(\mathbb{R}))$.
In order to define the concept of weak solution, we introduce the space $E(t_0,T_1)$ of functions $\varphi \in C^0([t_0,T_1],L^2(\mathbb{R}^n))$ satisfying

\medskip

\begin{itemize}
\item[$(i)$] $\varphi(\cdot,x) \in C^1((t_0,T_1),\cc)$ for all $x \in \mathbb{R}^n$, 
\item[$(ii)$] $\varphi(t,\cdot) \in C^2(\mathbb{R}^n,\cc)$ for all $t \in (t_0,T_1)$,
\item[$(iii)$] the functions $\partial_t \varphi+\langle B_0(t) x , \nabla_x \varphi \rangle$, $\nabla_x^2 \varphi$, $\varphi$ belong to  $L^2((t_0,T_1)\times\mathbb{R}^n)$.
\end{itemize}

\medskip

We consider the following notion of weak solution:

\medskip

\begin{definition}
Let $T_0 \leqslant t_0 < T_1$, $A_0, B_0 \in C^0([T_0,T_1],M_n(\mathbb{R}))$ and $k_0 \in L^2(\mathbb{R}^n)$.
A weak solution to the Cauchy problem \emph{(\ref{CY_Hom})} is a function $k \in C^0([t_0,T_1],L^2(\mathbb{R}^n))$ such that
$k(t_0)=k_0$ in $L^2(\mathbb{R}^n)$, and satisfying for all $\varphi \in E(t_0,T_1)$, $t_* \in (t_0,T_1)$,
$$\int\limits_{\mathbb{R}^n} \big(k(t_*,x) \varphi(t_*,x)-k_0(x)\varphi(t_0,x)\big) dx
= \int\limits_{t_0}^{t_*} \int\limits_{\mathbb{R}^n} k(t,x) \big(\partial_t \varphi(t,x)+P_0(t)^*\varphi(t,x) \big) dx dt,$$
with 
\begin{equation} \label{def:P*_appendix}
P_0(t)^*=\frac{1}{2}\emph{\text{Tr}}\big(A_0(t)A_0(t)^T\nabla_x^2\big)+\big\langle B_0(t)x,\nabla_x\big\rangle+\frac{1}{2}\emph{\text{Tr}}\big(B_0(t)\big), \quad t \in (T_0,T_1).
\end{equation}
\end{definition}

\medskip

We establish the following result:

\medskip

\begin{proposition} \label{Prop:WP_Hom}
Let $T_0<T_1$, $A_0, B_0 \in C^0([T_0,T_1],M_n(\mathbb{R}))$ and $\mathscr{T}=\{ (t,t_0) : T_0 \leqslant t_0 \leqslant t \leqslant T_1 \}$. 
There exists a strongly continuous mapping
$$\begin{array}{cccl}
U: & \mathscr{T} & \rightarrow & \mathcal{L}(L^2(\mathbb{R}^n)), \\
   & (t,t_0)     & \mapsto     & U(t,t_0),
\end{array}$$
with $\mathcal{L}(L^2(\mathbb{R}^n))$ denoting the space of bounded linear operators on $L^2(\mathbb{R}^n)$, satisfying

\medskip

\begin{itemize}
\item[$(i)$] $\forall T_0 \leqslant t \leqslant T_1, \quad U(t_0,t_0)=I$, 
\item[$(ii)$] $\forall T_0 \leqslant t_0 \leqslant t_1 \leqslant t_2 \leqslant T_1, \quad U(t_2,t_1) U(t_1,t_0) = U(t_2,t_0)$,
\item[$(iii)$] for all $k_0 \in L^2(\mathbb{R}^n)$, the function $k(t)=U(t,t_0)k_0$ is the unique weak solution to the Cauchy problem \emph{(\ref{CY_Hom})}.
\end{itemize}

\medskip

\noindent
Furthermore, the Fourier transform of the function $k(t)=U(t,t_0)k_0$ is given by
\begin{equation}\label{k_explicit}
\widehat{k}(t,\xi) = \widehat{k}_0\big(R_0(t,t_0)^T \xi \big)  e^{\frac{1}{2} \int_{t_0}^t\emph{\text{Tr}}(B_0(s))ds}e^{-\frac{1}{2} \int_{t_0}^t |A_0(s)^TR_0(t,s)^T \xi|^2ds},
\end{equation}
where $R_0$ denotes the resolvent associated to the linear time-varying system $\dot{X}(t)=B_0(t) X(t)$, that is, for all $T_0 \leqslant t_0\,, t \leqslant T_1$,
\begin{equation}\label{rt9}
\left\lbrace \begin{array}{ll}
\frac{\partial R_0}{\partial t} (t,t_0)=B_0(t) R_0(t,t_0), \\
R_0(t_0,t_0)=I_n.
\end{array} \right.
\end{equation}
\end{proposition}

\medskip

In the above statement, the normalization of the Fourier transform with respect to the space variable is given by
$$\widehat{k}(t,\xi)=\int_{\mathbb{R}^n} k(t,x) e^{-i x \cdot \xi} dx.$$
Following~\cite{pazy} (Chapter 5, Section~5.1, Definition~5.3, p.~129), the two parameter family of bounded linear operators $(U(t_1,t_2))_{(t_1,t_2) \in \mathscr{T}}$ is called the evolution system associated to the homogeneous equation (\ref{CY_Hom}). More specifically, we shall say that the mapping $U(t,t_0)$ is the evolution mapping associated to the family of operators $s \in [t_0,t] \mapsto P_0(s)$.

\begin{proof} Let $t_0 \in [T_0,T_1]$ and $k_0 \in L^2(\mathbb{R}^n)$.

\medskip

\noindent
\textit{Step 1.} We first derive heuristically an explicit expression of the Fourier transform $\widehat{k}$.
To that end, we consider $k$ a smooth solution to the Cauchy problem (\ref{CY_Hom}) and define the function $K : [t_0,T_1] \times \mathbb{R}^n \rightarrow \mathbb{C}$ by
\begin{equation} \label{CVAR:k/K}
k(t,x)=K\big(t, R_0(t_0,t)x\big).
\end{equation}
We recall for instance from~\cite{coron_book} (Proposition~1.5) that
\begin{equation}\label{rt1}
\forall T_0 \leq t_1,t_2,t_3 \leq T_1, \quad R_0(t_1,t_2)R_0(t_2,t_1)=I_n, \quad R_0(t_1,t_2)R_0(t_2,t_3)=R_0(t_1,t_3)
\end{equation}
and
\begin{equation}\label{rt2}
\forall T_0 \leq t_1,t_2 \leq T_1, \quad (\partial_2 R_0)(t_1,t_2)= - R_0(t_1,t_2)B_0(t_2).
\end{equation}
According to (\ref{rt1}) and (\ref{rt2}), the function $K$ is well-defined and a direct computation provides that 
\begin{equation}\label{rt3}
(\partial_t k)(t,x) + \big\langle B_0(t)x,(\nabla_x k)(t,x)\big\rangle = (\partial_t K)\big(t,R_0(t_0,t)x\big).
\end{equation}
It follows from (\ref{CY_Hom}) and (\ref{rt3}) that 
$$\left\lbrace \begin{array}{l}
\partial_t K(t,y) - \frac{1}{2} \text{Tr}\big(R_0(t_0,t) A_0(t) A_0(t)^T R_0(t_0,t)^T  \nabla_y^2 K(t,y) \big) + \frac{1}{2} \text{Tr}\big(B_0(t)\big) K(t,y)=0, \\
K|_{t=t_0} = k_0 \in L^2(\mathbb{R}^n).
\end{array} \right.$$
By taking the Fourier transform, we deduce that 
$$\left\lbrace \begin{array}{l}
\partial_t \widehat{K}(t,\eta) + \frac{1}{2}|A_0(t)^TR_0(t_0,t)^T \eta|^2 \widehat{K}(t,\eta)  + \frac{1}{2} \text{Tr}\big(B_0(t)\big) \widehat{K}(t,\eta)=0,  \\
\widehat{K}(t_0,\eta) = \widehat{k}_0(\eta).
\end{array} \right.$$
It leads to the following explicit expression
\begin{equation} \label{defK}
\forall (t,\eta) \in [t_0,T_1]\times\mathbb{R}^n, \quad \widehat{K}(t,\eta)=\widehat{k}_0(\eta) e^{- \frac{1}{2} \int_{t_0}^t |A_0(s)^T R_0(t_0,s)^T \eta |^2ds} e^{- \frac{1}{2} \int_{t_0}^t \text{Tr}(B_0(s))ds}.
\end{equation}
By using the Liouville formula 
\begin{equation}\label{iop1}
\forall t_1, t_2 \in [T_0,T_1], \quad \text{det}\big(R_0(t_2,t_1)\big)=\exp\Big(\int_{t_1}^{t_2}\textrm{Tr}\big(B_0(s)\big)ds\Big),
\end{equation}
see e.g.~\cite{cartan} (Proposition II.2.3.1), and the change of variable $y=R_0(t_0,t)x$, it follows that
\begin{multline*}
\widehat{k}(t,\xi) 
 =  \int_{\mathbb{R}^n} K\big(t,R_0(t_0,t)x\big) e^{-i x \cdot \xi} dx =  |\text{det}(R_0(t,t_0))| \int_{\mathbb{R}^n} K(t,y) e^{-i (R_0(t,t_0)y) \cdot \xi} dy \\
 =e^{\int_{t_0}^t \text{Tr}(B_0(s))ds} \widehat{K}\big(t,R_0(t,t_0)^T \xi \big) 
 = \widehat{k}_0 \big(R_0(t,t_0)^T \xi \big)  e^{\frac{1}{2} \int_{t_0}^t\text{Tr}(B_0(s))ds } e^{-\frac{1}{2} \int_{t_0}^t |A_0(s)^T R_0(t,s)^T\xi|^2ds}, 
\end{multline*}
since $R_0(t_0,s)^T R_0(t,t_0)^T = \big(R_0(t,t_0)R_0(t_0,s)\big)^T=R_0(t,s)^T$.
It proves the formula (\ref{k_explicit}).

\bigskip

\noindent
\textit{Step 2.} We prove that the $L^2$-function $k$ whose Fourier transform is given by (\ref{k_explicit}), is a weak solution to the Cauchy problem (\ref{CY_Hom}).
We easily notice that $k(t_0)=k_0$ and $k \in C^0([t_0,T_1],L^2(\mathbb{R}^n))$. Then, we use the change of function
\begin{equation}\label{rt5}
\varphi(t,x)=\psi\big(t,R_0(t_0,t)x\big)\big|\text{det}\big(R_0(t_0,t)\big)\big|.
\end{equation}
According to (\ref{rt1}), the function $\psi$ is well-defined.
It follows from the Liouville formula (\ref{iop1}) that 
\begin{multline}\label{rt6}
(\partial_t \varphi)(t,x) + \big\langle B_0(t)x,(\nabla_x\varphi)(t,x)\big\rangle\\
=\big|\text{det}\big(R_0(t_0,t)\big)\big|\Big(\partial_t\psi\big(t,R_0(t_0,t)x\big)-\text{Tr}\big(B_0(t)\big)\psi\big(t,R_0(t_0,t)x\big)\Big),
\end{multline}
since $\text{det}\big(R_0(t_0,t)\big) \in \rr_+^*$ for all $T_0 \leq t \leq T_1$.
According to (\ref{CVAR:k/K}), (\ref{rt5}) and (\ref{rt6}), it is sufficient to prove that for all $\psi \in \widetilde{E}(t_0,T_1)$, $t_* \in (t_0,T_1)$,
\begin{multline*}
\int_{\mathbb{R}^n} \big(K(t_*,y)\psi(t_*,y) - k_0(y)\psi(t_0,y)\big) dy \\
=  \int_{t_0}^{t_*} \int_{\mathbb{R}^n} K(t,y) \Big( \partial_t \psi 
+ \frac{1}{2} \text{Tr}\big(R_0(t_0,t)A_0(t)A_0(t)^TR_0(t_0,t)^T\nabla_y^2\psi\big)-\frac{1}{2}\text{Tr}\big(B_0(t)\big)\psi\Big)(t,y)dydt.
\end{multline*}
where $\widetilde{E}(t_0,T_1)$ stands for the space of functions $\psi \in C^0([t_0,T_1],L^2(\mathbb{R}^n))$ satisfying

\medskip

\begin{itemize}
\item[$(i)$] $\psi(\cdot,y) \in C^1((t_0,T_1),\cc)$ for all $y \in \mathbb{R}^n$,
\item[$(ii)$] $\psi(t,\cdot) \in C^2(\mathbb{R}^n,\cc)$ for all $t \in (t_0,T_1)$,
\item[$(iii)$] the functions $\partial_t \psi$, $\nabla_y^2 \psi$, $\psi$ belong to $L^2((t_0,T_1)\times\mathbb{R}^n)$.
\end{itemize}

\medskip

\noindent
For all $\psi \in \widetilde{E}(t_0,T_1)$ and $t^* \in (t_0,T_1)$, it follows from the Plancherel theorem, (\ref{k_explicit}) and (\ref{CVAR:k/K}) that
\begin{align*}
  & \int_{t_0}^{t_*} \int_{\mathbb{R}^n} K(t,y) 
    \Big( \partial_t \psi + \frac{1}{2} \text{Tr}\big(R_0(t_0,t)A_0(t)A_0(t)^TR_0(t_0,t)^T \nabla_y^2\psi\big) - \frac{1}{2}\text{Tr}\big(B_0(t)\big) \psi \Big)(t,y) dy dt  \\
= & \frac{1}{(2\pi)^n} \int\limits_{t_0}^{t_*} \int\limits_{\mathbb{R}^n} \widehat{K}(t,\eta) 
    \Big(\overline{ \partial_t \widehat{\overline{\psi}} - \frac{1}{2} |A_0(t)^TR_0(t_0,t)^T  \eta|^2 \widehat{\overline{\psi}}-\frac{1}{2}\text{Tr}\big(B_0(t)\big) \widehat{\overline{\psi}} }\Big)(t,\eta)  d\eta dt  \\
= & \frac{1}{(2\pi)^n}  \int\limits_{t_0}^{t_*} \int\limits_{\mathbb{R}^n} \frac{\partial}{\partial t} \Big[ \widehat{K}(t,\eta) \overline{\widehat{\overline{\psi}}(t,\eta)} \Big] d\eta dt 
=  \frac{1}{(2\pi)^n} \int\limits_{\mathbb{R}^n} \Big( \widehat{K}(t_*,\eta) \overline{\widehat{\overline{\psi}}(t_*,\eta)} - \widehat{K}(t_0,\eta) \overline{\widehat{\overline{\psi}}(t_0,\eta)} \Big) d\eta \\
= & \int\limits_{\mathbb{R}^n} \big( K(t_*,y) \psi(t_*,y) - K(t_0,y) \psi(t_0,y) \big) dy=\int\limits_{\mathbb{R}^n} \big( K(t_*,y) \psi(t_*,y) - k_0(y) \psi(t_0,y) \big) dy.
\end{align*}

\bigskip

\noindent
\textit{Step 3: Definition and properties of the evolution system.}
For all $(t,t_0) \in \mathscr{T}$ and $k_0 \in L^2(\mathbb{R}^n)$, we define $U(t,t_0)k_0$ as the $L^2$-function $k(t)$ whose Fourier transform is given by (\ref{k_explicit}). With this definition, we easily check that $U(t_0,t_0)=I$ for all $T_0 \leq t_0 \leq T_1$ and that the mapping $U$ is strongly continuous from $\mathscr{T}$ to $\mathcal{L}(L^2(\mathbb{R}^n))$ thanks to Plancherel theorem. On the other hand, with $k_1=U(t_1,t_0)k_0$, $k_2=U(t_2,t_0)k_0$ and $k_3=U(t_2,t_1)k_1$, it follows from (\ref{k_explicit}) that for all $T_0 \leqslant t_0 \leqslant t_1 \leqslant t_2 \leqslant T_1$, $k_0 \in L^2(\mathbb{R}^n)$, 
$$\widehat{k}_1(\xi) = \widehat{k}_0\big(R_0(t_1,t_0)^T\xi\big)e^{\frac{1}{2}\int_{t_0}^{t_1}\text{Tr}(B_0(s))ds}e^{-\frac{1}{2}\int_{t_0}^{t_1}|A_0(s)^TR_0(t_1,s)^T \xi|^2ds},$$
$$\widehat{k}_2(\xi) = \widehat{k}_0\big(R_0(t_2,t_0)^T\xi\big)e^{\frac{1}{2}\int_{t_0}^{t_2}\text{Tr}(B_0(s))ds}e^{-\frac{1}{2}\int_{t_0}^{t_2}|A_0(s)^TR_0(t_2,s)^T\xi|^2ds}$$
and
\begin{align*}
\widehat{k}_3(\xi)=& \ \widehat{k}_1\big(R_0(t_2,t_1)^T\xi\big)e^{\frac{1}{2} \int_{t_1}^{t_2}\text{Tr}(B_0(s))ds}e^{-\frac{1}{2} \int_{t_1}^{t_2}|A_0(s)^TR_0(t_2,s)^T \xi|^2ds} \\
= & \  \widehat{k}_0\big(R_0(t_1,t_0)^TR_0(t_2,t_1)^T\xi\big)
 e^{\frac{1}{2} \int_{t_0}^{t_2}\text{Tr}(B_0(s))ds}\\
 & \ \times e^{-\frac{1}{2} \int_{t_0}^{t_1}|A_0(s)^TR_0(t_1,s)^T R_0(t_2,t_1)^T\xi|^2ds} 
 e^{-\frac{1}{2} \int_{t_1}^{t_2}|A_0(s)^TR_0(t_2,s)^T\xi|^2ds}\\
 = & \ \widehat{k}_2(\xi),
\end{align*}
since 
$$R_0(t_1,s)^TR_0(t_2,t_1)^T=\big(R_0(t_2,t_1)R_0(t_1,s)\big)^T=R_0(t_2,s)^T$$
and
$$R_0(t_1,t_0)^TR_0(t_2,t_1)^T=\big(R_0(t_2,t_1)R_0(t_1,t_0)\big)^T=R_0(t_2,t_0)^T.$$
It proves that for all $T_0 \leqslant t_0 \leqslant t_1 \leqslant t_2 \leqslant T_1$,
$$U(t_2,t_1)U(t_1,t_0)=U(t_2,t_0).$$

\bigskip

\noindent
\textit{Step 4: Uniqueness of the weak solution to the Cauchy problem \emph{(\ref{CY_Hom})}.}
Let $k$ be a  weak solution to the Cauchy problem (\ref{CY_Hom}) associated with the initial datum $k_0=0$. It follows that for all $\varphi \in E(t_0,T_1)$, $t_* \in (t_0,T_1)$,
\begin{equation} \label{WS_uniqueness}
\int\limits_{\mathbb{R}^n} k(t_*,x) \varphi(t_*,x) dx
= \int\limits_{t_0}^{t_*} \int\limits_{\mathbb{R}^n} k(t,x) \big( \partial_t \varphi(t,x) + P_0(t)^* \varphi(t,x) \big) dx dt.
\end{equation}
Let $t_* \in [t_0,T_1]$ be fixed. We aim at proving that $k(t_*)= 0$. To that end,
we consider a sequence $(g_p)_{p \geq 1}$ of $C_0^{\infty}(\rr^n)$ functions satisfying 
$$\lim_{p \to +\infty}\|\widehat{g}_p - \widehat{k}(t_*)\|_{L^2(\mathbb{R}^n)}.$$
By Plancherel theorem, we observe that 
\begin{equation}\label{rt7}
\lim_{p \to +\infty}\|g_p - k(t_*)\|_{L^2(\mathbb{R}^n)}=0.
\end{equation}
Following the very same strategy as in the two first steps, we build a weak solution $\varphi_p:(t_0,t_*) \times \mathbb{R}^n \rightarrow \mathbb{C}$ to the Cauchy problem
$$\left\lbrace \begin{array}{ll}
\partial_t \varphi_p(t,x) + P_0(t)^* \varphi_p(t,x)=0, \\
\varphi_p|_{t=t_*}=\overline{g_p}.
\end{array}\right.$$
By deriving a similar formula as in (\ref{k_explicit}), we notice that the function $\varphi_p$ is smooth in the space variable as its Fourier transform in the space variable is compactly supported. This similar formula as in (\ref{k_explicit}) also shows that the function $\varphi_p$ is smooth in the time variable. It follows that the function $\varphi_p$ is a pointwise solution of the equation:
$$\forall (t,x) \in (t_0,t^*)\times\mathbb{R}^n, \quad \partial_t \varphi_p(t,x) + P_0(t)^* \varphi_p(t,x)=0.$$
Furthermore, we easily check that $\varphi_p$ is an admissible test function.
Then, we deduce from (\ref{WS_uniqueness}) and (\ref{rt7}) that
$$\forall p \geq 1, \quad \int\limits_{\mathbb{R}^n} k(t_*,x)\overline{g_p(x)}dx=(k(t_*),g_p)_{L^2(\rr^n)}=0,$$
implying that $k(t_*)=0$ when passing to the limit $p \to +\infty$. It ends the proof of the uniqueness of the weak solution to the Cauchy problem (\ref{CY_Hom}).
\end{proof}

\bigskip

Regarding the non-homogeneous equation
\begin{equation} \label{CY_Nonhom}
\left\lbrace \begin{array}{ll}
\partial_t h(t,x) - P_0(t) h(t,x)=u(t,x), \quad & (t,x) \in (t_0,T_1) \times \mathbb{R}^n, \\
h|_{t=t_0} = h_0 \in L^2(\mathbb{R}^n),
\end{array} \right.
\end{equation}
we use the notion of mild solutions defined in~\cite{pazy} (Chapter~5, Section~5.5, Definition~5.1, p.~146):

\medskip

\begin{definition}
Let $T_0 \leqslant t_0 < T_1$, $A_0, B_0 \in C^0([T_0,T_1],M_n(\mathbb{R}))$, $k_0 \in L^2(\mathbb{R}^n)$ and $u \in L^1((t_0,T_1),L^2(\mathbb{R}^n))$.
The mild solution to the Cauchy problem \emph{(\ref{CY_Nonhom})} is the function $h \in C^0([t_0,T_1],L^2(\mathbb{R}^n))$  given by
$$h(t)=U(t,t_0)h_0 + \int_{t_0}^t U(t,s)u(s)ds,$$
with equality in $L^2(\mathbb{R}^n)$ for all $t \in [t_0,T_1]$, where $U$ stands for the evolution system given by 
Proposition~\ref{Prop:WP_Hom}.
\end{definition}

\subsection{Explicit computation of a solution to a particular Cauchy problem}
Let $A, B \in C^0(I,M_{n}(\rr))$ with $I$ being an open interval of $\rr$ containing zero and $T>0$. We deduce from Proposition~\ref{Prop:WP_Hom} that the Fourier transform of the mild solution to the Cauchy problem on $(t_0,T) \times \rr^n$,
$$\left\lbrace \begin{array}{l}
\partial_t g(t,x) - \frac{1}{2}\textrm{Tr}\big(A(T-t)A(T-t)^T\nabla_x^2 g(t,x)\big) + \big\langle B(T-t)x , \nabla_x g(t,x) \big\rangle \\
\hspace{9cm}+\frac{1}{2}\textrm{Tr}\big(B(T-t)\big) g(t,x)= 0\,,\\
g(t_0,x)=g_0(x),
\end{array}\right.$$
where $g_0 \in L^2(\mathbb{R}^n)$, with $t_0 \in [0,T]$, is given by
\begin{equation}\label{sdf1}
\widehat{g}(t,\xi)=\exp\Big(\frac{1}{2}\int_{t_0}^t\textrm{Tr}\big(B(T-s)\big)ds\Big)
 \widehat{g}_0(R(t,t_0)^T\xi) e^{-\frac{1}{2}\int_{t_0}^t |A(T-s)^TR(t,s)^T\xi|^2 ds},
\end{equation}
where the mapping
$$\begin{array}{ll}
R : [0,T] \times [0,T] & \rightarrow M_n(\rr), \\
\ \ \ \  (t_1,t_2)  & \mapsto R(t_1,t_2),
\end{array}$$
stands for the resolvent of the time-varying linear system 
$$\overset{.}{X}(t)=B(T-t)X(t).$$

\subsection{Hilbert uniqueness method}

This section is devoted to the proof of the following result:

\medskip

\begin{proposition}\label{Prop:HUM}
Let $T_0<0<T_1$, $A_0, B_0 \in C^0([T_0,T_1],M_n(\mathbb{R}))$, 
$P_0(t)$ be the non-autonomous Ornstein-Uhlenbeck operator defined in (\ref{def:P_appendix}) and $P_0(t)^*$ its adjoint given in (\ref{def:P*_appendix}).
The null-controllability from the set $\omega$ in time $0<T<T_1$ of the system
\begin{equation} \label{syst_controle}
\left\lbrace \begin{array}{ll}
\partial_t f(t,x) - P_0(t) f(t,x)=u(t,x)\un_{\omega}(x), \quad & (t,x) \in (0,T_1) \times \mathbb{R}^n, \\
f|_{t=0} = f_0 \in L^2(\mathbb{R}^n),
\end{array} \right.
\end{equation}
is equivalent to the existence of an observability constant $C>0$ such that,
for all $g_0 \in L^2(\mathbb{R}^n)$, the mild solution to the Cauchy problem
\begin{equation} \label{syst_adj}
\left\lbrace \begin{array}{ll}
\partial_t g(t,x) - P_0(T-t)^* g(t,x)=0, \quad & (t,x) \in (0,T-T_0) \times \mathbb{R}^n, \\
g|_{t=0} = g_0 \in L^2(\mathbb{R}^n),
\end{array} \right.
\end{equation}
satisfies
$$\int\limits_{\mathbb{R}^n} |g(T,x)|^2 dx \leqslant C \int\limits_0^T \int\limits_\omega |g(t,x)|^2 dx dt.$$
\end{proposition}

\medskip

\noindent
Instrumental in the proof of Proposition~\ref{Prop:HUM} is the following result:

\medskip

\begin{lemma} \label{Prop:U*}
If $U$ stands for the evolution system given by Proposition \ref{Prop:WP_Hom}, then the $L^2$-adjoint $U(t,t_0)^*$ of the evolution mapping $U(t,t_0)$ is  equal to the evolution mapping $\tilde{U}(t-t_0,0)$ associated to the family of operators $s \in [0,t-t_0] \mapsto P_0(t-s)^*$.
\end{lemma}

\medskip

\begin{proof} Let $T_0 \leqslant t_0 < t \leqslant T_1$ and $g_0 \in L^2(\mathbb{R}^n)$. Setting $g(t-t_0)=\tilde{U}(t-t_0,0)g_0$, we deduce from 
Proposition~\ref{Prop:WP_Hom} with the suitable substitutions that 
\begin{multline}\label{rt11}
\widehat{g}(t-t_0,\xi) =\\
  \widehat{g}_0 \big( \mathcal{R}(t-t_0,0)^T \xi \big)  
e^{-\frac{1}{2} \int_{0}^{t-t_0}\text{Tr}(B_0(t-s))ds}e^{-\frac{1}{2} \int_{0}^{t-t_0}|A_0(t-s)^T \mathcal{R}(t-t_0,s)^T \xi|^2ds},
\end{multline}
where $\mathcal{R}$ stands for the resolvent associated to the system $\dot{X}(s)=-B_0(t-s)X(s)$, that is,
$$\left\lbrace \begin{array}{l}
\frac{\partial \mathcal{R}}{\partial s}(s,s_0)=-B_0(t-s)\mathcal{R}(s,s_0), \\
\mathcal{R}(s_0,s_0)=I_n.
\end{array}\right.$$
According to (\ref{rt9}), we notice that
\begin{equation}\label{rt10}
\mathcal{R}(s,s_0)=R_0(t-s,t-s_0).
\end{equation}
It follows from (\ref{rt11}) and (\ref{rt10}) that 
\begin{equation}\label{rt13}
\widehat{g}(t-t_0,\xi)  = \widehat{g_0} \big(R_0(t_0,t)^T \xi \big) 
e^{-\frac{1}{2} \int_{t_0}^{t}\text{Tr}(B_0(s))ds}e^{-\frac{1}{2} \int_{t_0}^{t}|A_0(s)^TR_0(t_0,s)^T \xi|^2ds}.
\end{equation}
We deduce from Plancherel theorem, (\ref{k_explicit}), the change of variable $\eta= R_0(t,t_0)^T \xi$,
the Liouville formula (\ref{iop1}) and (\ref{rt13}) that for all $k_0$, $g_0 \in L^2(\mathbb{R}^n)$,
\begin{align*}
  & \big(U(t,t_0)k_0,g_0\big)_{L^2(\mathbb{R}^n)}
= \frac{1}{(2\pi)^n}\int_{\mathbb{R}^n} \widehat{k}(t,\xi) \overline{\widehat{g}_0(\xi)} d\xi 
\\
= & \frac{1}{(2\pi)^n}\int_{\mathbb{R}^n} \widehat{k}_0 \big(R_0(t,t_0)^T \xi \big)  
  e^{\frac{1}{2} \int_{t_0}^t\text{Tr}(B_0(s))ds}
  e^{-\frac{1}{2} \int_{t_0}^t|A_0(s)^TR_0(t,s)^T\xi|^2ds} \overline{\widehat{g}_0(\xi)} d\xi 
\\
= & \frac{1}{(2\pi)^n}\int_{\mathbb{R}^n} \widehat{k}_0(\eta)   
e^{-\frac{1}{2} \int_{t_0}^t\text{Tr}(B_0(s))ds}
  e^{-\frac{1}{2} \int_{t_0}^t|A_0(s)^TR_0(t,s)^TR_0(t_0,t)^T\eta|^2ds} 
 \overline{\widehat{g}_0\big( R(t_0,t)^T \eta\big)} d\eta, 
\\
= & \frac{1}{(2\pi)^n}\int_{\mathbb{R}^n} \widehat{k}_0(\eta)\overline{\widehat{g}(t-t_0,\eta)} d\eta=\int_{\mathbb{R}^n} k_0(x)\overline{g(t-t_0,x)}dx
=\big(k_0,\tilde{U}(t-t_0,0)g_0\big)_{L^2(\rr^n)},
\end{align*}
since
$$R_0(t,s)^TR_0(t_0,t)^T=R_0(t_0,s)^T, \quad \big|\text{det}\big(R_0(t_0,t)\big)\big|=\exp\Big(-\int_{t_0}^{t}\textrm{Tr}\big(B_0(s)\big)ds\Big).$$
This ends the proof of Lemma~\ref{Prop:U*}.
\end{proof}

\begin{proof} We can now resume the proof of Proposition~\ref{Prop:HUM}. We consider the following linear mappings
$$
\begin{array}{crcl}
C_2: & L^2(\mathbb{R}^n) & \rightarrow & L^2(\mathbb{R}^n), \\
     &     f_0           & \mapsto     & U(T,0)f_0,
\end{array}$$
and
$$\begin{array}{crcl}
C_3: & L^2((0,T)\times\omega) & \rightarrow & L^2(\mathbb{R}^n), \\
     &     u                  & \mapsto     & \int_0^T U(T,s) u(s) ds.
\end{array}$$
For any $f_0 \in L^2(\mathbb{R}^n)$, its image by the first mapping $C_2(f_0)=k(T)$ is the value at time $T$ of the weak solution to the Cauchy problem (\ref{CY_Hom}) with $t_0=0$ associated to the initial datum $k_0=f_0$. On the other hand, for any $u \in L^2((0,T)\times\omega)$, its image by the second mapping $C_3(u)=h(T)$ is the value at time $T$ of the mild solution to the Cauchy problem (\ref{CY_Nonhom}) with $t_0=0$ and $h_0=0$.

The null-controllability from the set $\omega$ in time $T$ of the system (\ref{syst_controle}) is equivalent to the inclusion 
$$C_2(L^2(\mathbb{R}^n)) \subset C_3(L^2((0,T)\times\omega)),$$ 
since $f(T)=C_2(f_0)+C_3(u)$.
According to~\cite{coron_book} (Lemma~2.48), this is also equivalent to the existence of a positive constant $M>0$ such that for all $g_0 \in L^2(\mathbb{R}^n)$,
\begin{equation}\label{rt20}
\| C_2^*g_0 \|_{L^2(\mathbb{R}^n)} \leqslant M \|C_3^* g_0 \|_{L^2((0,T)\times\omega)}.
\end{equation}
We deduce from Lemma~\ref{Prop:U*} that 
\begin{equation}\label{rt21}
C_2^* g_0 = U(T,0)^* g_0= \tilde{U}(T,0)g_0 = g(T),
\end{equation}
where $g$ is the weak solution of (\ref{syst_adj}). On the other hand, we have for all $u \in L^2((0,T)\times\omega)$,
\begin{multline}\label{rt22}
(C_3 u,g_0)_{L^2(\mathbb{R}^n)}
= \int_0^T \big(U(T,s)u(s),g_0\big)_{L^2(\mathbb{R}^n)} ds
= \int_0^T \big(u(s),U(T,s)^*g_0\big)_{L^2(\omega)} ds\\
= \int_0^T \big(u(s),\tilde{U}(T-s,0)g_0\big)_{L^2(\omega)} ds= \int_0^T \big(u(s),g(T-s)\big)_{L^2(\omega)} ds,
\end{multline}
where $g$ is the weak solution of (\ref{syst_adj}). It follows from (\ref{rt21}) and (\ref{rt22}) that the estimate (\ref{rt20}) reads as 
$$\int\limits_{\mathbb{R}^n} |g(T,x)|^2 dx \leqslant M \int\limits_0^T \int\limits_\omega |g(t,x)|^2 dx dt.$$
This ends the proof of Proposition~\ref{Prop:HUM}.
\end{proof}

\subsection{Counter-example to the necessity of the Kalman type condition for null-controllability}\label{counter}
The Kalman type condition (\ref{wer1}) is a sufficient condition for the null-controllability of 
the non-autonomous Ornstein-Uhlenbeck equation (\ref{syst_LR1}). As shown by the following example in dimension $n=2$, this condition is in general not necessary to get null-controllability.

Let $T>0$ and $B=0$. The resolvent defined in (\ref{resol}) is given by $R(t_1,t_2)=I_2$, for all $t_1,t_2 \in [0,T]$.
We consider smooth real-valued functions $f, g \in C_0^\infty((0,1),\mathbb{R})$ satisfying $f = 0$ on $[1/2,1]$,
$g = 0$ on $[0,1/2]$ and 
$$\int_0^1 f(t)^2 dt = \int_0^1 g(t)^2 dt = 1.$$
We define
$$\mathcal{A}(t)=\left( \begin{array}{cc}
f(t) & 0 \\ 0 & g(t)
\end{array} \right)\,, \quad 0 \leq t \leq 1.$$
We observe that the matrices $\frac{d^k}{dt^k} \mathcal{A}(t)$, with $k \geq 0$, are either zero or rank one for all $0 \leq t \leq 1$.
The time-dependent matrix $A$ defining the Ornstein-Uhlenbeck operator is devised piecewise according to the partition of the interval $[0,T]$ used in the proof of Theorem~\ref{meta_thm}.

Let $0 <\rho < 1/3$ and $K>0$ be the positive constant verifying
$$\sum\limits_{k=1}^{+\infty} \frac{2K}{4^{k\rho}} = T\,.$$
We define for all $k \geqslant 1$,
$$\tau_k=\frac{K}{4^{k\rho}}\,, \quad \alpha_0=0, \quad
\alpha_k=\sum_{j=1}^{k} 2 \tau_j\,,$$ 
$$J_k=\left[T - \alpha_{k-1}-\tau_k , T - \alpha_{k-1} \right]\,, \quad 
I_k=\left[T-\alpha_k , T - \alpha_k + \tau_k \right],$$
so that $[T-\alpha_k ,T-\alpha_{k-1}]=I_k \cup J_k$.
Let $A$ be defined on $[0,T]$ as
$$A(T-s)= \left\lbrace \begin{array}{ll}
\tau_k \mathcal{A} \left( \frac{s-(T-\alpha_k)}{\tau_k} \right) \quad & \text{ if } s \in I_k\,, \\
0                                                                     & \text{ if } s \in J_k\,. \\
\end{array}\right.$$
The Kalman type condition (\ref{wer1}) does not hold as
$$\forall t_0 \in [0,T], \quad \textrm{dim}(\textrm{Span}\{\tilde{A}_{k}(t_0)x : \ x \in \rr^2,\ k \geq 0\}) \leq 1,$$
where $\widetilde{A}_k(t)=\frac{d^k}{dt^k} A(t)$, if $0 \leq t \leq T$ and $k \geq 0$.
On the other hand, we observe that for all $0 \leq \tau \leq T$, $\xi=(\xi_1,\xi_2) \in \rr^2$,
$$\begin{array}{ll}
  & \int\limits_{I_k} |A(T-s)^T R(\tau,s)^T \xi |^2 ds  \\ 
= &  \xi_1^2 \, \int\limits_{T-\alpha_k}^{T-\alpha_k+\frac{\tau_k}{2}} \left| \tau_k f\left( \frac{s-(T-\alpha_k)}{\tau_k} \right) \right|^2 ds \,
    + \, \xi_2^2 \, \int\limits_{T-\alpha_k+\frac{\tau_k}{2}}^{T-\alpha_k+\tau_k} \left| \tau_k g\left( \frac{s-(T-\alpha_k)}{\tau_k} \right) \right|^2 ds  \\
= & \tau_k^3 |\xi|^2\,.
\end{array}$$
By taking $n=2$, $\eps=0$, $\tilde{\eps}=T$ and substituting $\sum_{j=1}^k\tau_j^3$ to the term $\alpha_k^{2n-1}$,
the proof of Theorem~\ref{meta_thm} can be readily adapted 
in order to prove that the non-autonomous Ornstein-Uhlenbeck equation posed in the $L^2(\rr^2)$-space
$$\left\lbrace \begin{array}{l}
\partial_t f(t,x) - \frac{1}{2}\textrm{Tr}\big(A(t)A(t)^T\nabla_x^2 f(t,x)\big)=u(t,x)\un_{\omega}(x), \\
f|_{t=0}=f_0 \in L^2(\rr_x^2),                                       
\end{array}\right.$$
is null-controllable from the set $\omega$ in time greater or equal to $T$, if $\omega$ is an open subset of $\rr^2$ satisfying (\ref{hyp_omega}).

\end{document}